\documentclass{amsart}
\usepackage{amssymb}
\usepackage{amsmath}
\usepackage{latexsym}
\usepackage{eepic}
\usepackage{epsfig}
\usepackage{graphicx}
\usepackage{pb-diagram,lamsarrow,pb-lams,amscd}
\usepackage{eufrak}
\usepackage{mathrsfs}
\usepackage[english]{babel}
\usepackage[all,cmtip]{xy}

\DeclareMathAlphabet{\mathpzc}{OT1}{pzc}{m}{it}
\newtheorem{theorem}{Theorem}[section]
\newtheorem{theorem-definition}[theorem]{Theorem-Definition}
\newtheorem{lemma-definition}[theorem]{Lemma-Definition}
\newtheorem{definition-prop}[theorem]{Proposition-Definition}

\newtheorem{prop}[theorem]{Proposition}
\newtheorem{lemma}[theorem]{Lemma}
\newtheorem{cor}[theorem]{Corollary}
\newtheorem{definition}[theorem]{Definition}

\newtheorem{conjecture}[theorem]{Conjecture}

\newtheorem{question}[theorem]{Question}
\newenvironment{remark}{\vspace{4pt}\noindent\textbf{Remark.}}{\qed\vspace{4pt}}

\newcommand{\LL}{\ensuremath{\mathbb{L}}}

\newcommand{\N}{\ensuremath{\mathbb{N}}}
\newcommand{\Z}{\ensuremath{\mathbb{Z}}}
\newcommand{\Q}{\ensuremath{\mathbb{Q}}}

\newcommand{\R}{\ensuremath{\mathbb{R}}}
\newcommand{\C}{\ensuremath{\mathbb{C}}}

\newcommand{\A}{\ensuremath{\mathbb{A}}}

\newcommand{\X}{\ensuremath{\mathscr{X}}}

\renewcommand{\R}{\ensuremath{\mathbb{R}}}
\renewcommand{\C}{\ensuremath{\mathbb{C}}}

\renewcommand{\A}{\ensuremath{\mathbb{A}}}
\renewcommand{\X}{\ensuremath{\mathfrak{X}}}

\newcommand{\mU}{\ensuremath{\mathfrak{U}}}

\newcommand{\Spec}{\ensuremath{\mathrm{Spec}\,}}

\numberwithin{equation}{section}
\begin{document}
\title[Motivic zeta functions of abelian varieties]{Motivic zeta functions of abelian varieties, and the monodromy conjecture}
\author{Lars Halvard Halle}
\address{Institut f{\"u}r Algebraische Geometrie\\
Gottfried Wilhelm Leibniz Universit{\"a}t Hannover\\ Welfengarten
1
\\
30167 Hannover \\ Deutschland} \email{halle@math.uni-hannover.de}
\author[Johannes Nicaise]{Johannes Nicaise}
\address{Universit\'e Lille 1\\
Laboratoire Painlev\'e, CNRS - UMR 8524\\ Cit\'e Scientifique\\59655 Villeneuve d'Ascq C\'edex \\
France} \email{johannes.nicaise@math.univ-lille1.fr}

\thanks{The first author was partially supported by the Fund for Scientific
Research-Flanders (G.0318.06) and by DFG under grant Hu 337/6-1.
 The second author was partially supported by ANR-06-BLAN-0183.
}

\begin{abstract}
We prove for abelian varieties a global form of Denef and Loeser's
motivic monodromy conjecture, in arbitrary characteristic. More
precisely, we prove that for every tamely ramified abelian variety
$A$ over a complete discretely valued field, its motivic zeta
function has a unique pole at Chai's base change conductor $c(A)$
of $A$, and that the order of this pole equals one plus the
potential toric rank of $A$. Moreover, we show that for every
embedding of $\Q_\ell$ in $\C$, the value $\exp(2\pi i c(A))$ is
an $\ell$-adic tame monodromy eigenvalue of $A$. The main tool in
the paper is Edixhoven's filtration on the special fiber of the
N\'eron model of $A$, which measures the behaviour of the N\'eron
model under tame base change.
\end{abstract}

 \maketitle
\section{Introduction}
Let $K$ be a complete discretely valued field with ring of
integers $R$ and algebraically closed residue field $k$. We denote
by $p$ the characteristic exponent of $k$, and by $\N'$ the set of
strictly positive integers that are prime to $p$. We fix a prime
$\ell\neq p$.

 Given a smooth, proper and connected $K$-variety $X$
and a gauge form $\omega$ on $X$ (i.e., a nowhere vanishing
differential form of maximal degree), one can associate to the
pair $(X,\omega)$ a \textit{motivic generating series}
$S(X,\omega;T)$ as follows. For each $d\in \N'$, the field $K$
admits an extension $K(d)$ of degree $d$ which is unique up to
$K$-isomorphism. If we put $X(d)=X\times_K K(d)$ and if we denote
by $\omega(d)$ the pull-back of $\omega$ to $X(d)$, then the
motivic generating series $S(X,\omega;T)$ is given by
$$S(X,\omega;T)=\sum_{d\in \N'}\left(\int_{X(d)}|\omega(d)|\right)
T^d\in \mathcal{M}_k[[T]]$$ Here $\mathcal{M}_k$ denotes the
localized Grothendieck ring of $k$-varieties, and
$$\int_{X(d)}|\omega(d)|\in \mathcal{M}_k$$ is the motivic integral of the gauge
form $\omega(d)$ on $X(d)$. These motivic integrals were defined
in \cite{motrigid} and can be computed on a weak N\'eron model of
$X(d)$ (see Proposition \ref{prop-motint}).

The aim of this paper is to study the series $S(A,\omega;T)$ when
$A$ is a tamely ramified abelian $K$-variety and $\omega$ is a
``normalized Haar measure'' on $A$, i.e., a translation-invariant
gauge form that extends to a relative gauge form on the N\'eron
model of $A$. Such a normalized Haar measure always exists, and it
is unique up to multiplication with a unit in $R$. It follows that
$S(A,\omega;T)$ only depends on $A$. We denote the series
$\LL^{\mathrm{dim}\,(A)}\cdot S(A,\omega;T)$ by $Z_A(T)$, and we
call it the \textit{motivic zeta function} of the abelian
$K$-variety $A$\footnote{This notion should not be confused with
Kapranov's motivic zeta function of $A$, which is studied, for
instance, in \cite{bitt-abelian}. The two are not related in any
direct way.}.

It is easily seen that
\begin{equation}\label{intro-zeta}
Z_A(T)=\sum_{d\in
\N'}[\mathcal{A}(d)_s]\LL^{-ord_{\mathcal{A}(d)_s^o}(\omega(d))}T^d\quad
\in \mathcal{M}_k[[T]] \end{equation} where $\mathcal{A}(d)$ is
the N\'eron model of $A(d)$, and
$ord_{\mathcal{A}(d)_s^o}(\omega(d))$ denotes the order of the
gauge form $\omega(d)$ on $A(d)$ along the identity component
$\mathcal{A}(d)_s^o$ of the special fiber of $\mathcal{A}(d)$. So
the zeta function $Z_A(T)$ depends on two (related) factors: the
behaviour of the N\'eron model of $A$ under tame base change, and
the function $d\mapsto ord_{\mathcal{A}(d)_s^o}(\omega(d))$.

 Our analysis of both
factors heavily relies upon Edixhoven's results in \cite{edix}. He
constructs a filtration on the special fiber $\mathcal{A}_s$ of
the N\'eron model of $A$, which measures the behaviour of the
N\'eron model under tame base change. We prove that
$ord_{\mathcal{A}(d)_s^o}(\omega(d))$ can be expressed in terms of
the jumps in this filtration. We also prove that the sum of the
jumps is equal to the base change conductor $c(A)$ introduced by
Chai in \cite{chai}.

In \cite{edix}, Edixhoven shows as well that the N\'eron model of
$A$ is canonically isomorphic to the $G(K(d)/K)$-invariant part of
the Weil restriction of the N\'eron model of $A(d)$. This result
allows us to analyze the behaviour of $[\mathcal{A}(d)_s]\in
\mathcal{M}_k$ when $d$ varies. If we denote by $\phi_A(d)$ the
number of connected components of $\mathcal{A}(d)_s$, then we have
 $[\mathcal{A}(d)_s]=\phi_A(d)[\mathcal{A}(d)_s^o]$. We show that $[\mathcal{A}(d)_s^o]$ only depends
on $d$ modulo $e$, with $e$ the degree of the minimal extension of
$K$ where $A$ acquires semi-abelian reduction. The behaviour of
$\phi_A(\cdot)$ is interesting in its own right; we studied it in
the independent paper \cite{HaNi-comp}.

The first question that arises is the rationality of $Z_A(T)$, and
the nature of its poles. We prove that $Z_A(T)$ is rational, and
belongs to the subring
$$\mathcal{M}_k\left[T,\frac{1}{1-\LL^{a}T^{b}}\right]_{(a,b)\in \Z\times \Z_{>0},\ a/b=c(A)}$$
of $\mathcal{M}_k[[T]]$. In particular, the series $Z_A(\LL^{-s})$
has a unique pole at $s=c(A)$. We show that the order of this pole
 equals one plus the potential toric rank of $A$. Our proof does not use
resolution of singularities.
%

A second question we consider is the relation between the pole
$s=c(A)$ and the monodromy eigenvalues of $A$. To explain the
motivation behind this question, we need some background. Let $X$
be a smooth connected complex variety, endowed with a dominant
morphism $f:X\rightarrow \Spec \C[t]$. Let $x$ be a closed point
of the special fiber of $f$. We denote by $\mathscr{F}_x$ the
analytic Milnor fiber of $f$ at $x$ \cite[\S\,9.1]{NiSe}. It is a
separated smooth rigid $\C((t))$-variety, and it serves as a
non-archimedean model for the topological Milnor fibration of $f$
at $x$. If $\phi/dt$ is a so-called \textit{Gelfand-Leray form} on
$\mathscr{F}_x$ \cite[9.5]{NiSe} then the series
$S(\mathscr{F}_x,\phi/dt;T)\in \mathcal{M}_{\C}[[T]]$ can be
defined in a similar way as above. Up to normalization, it
coincides with Denef and Loeser's motivic zeta function of $f$ at
$x$ \cite[9.7]{ni-trace}. Denef and Loeser's monodromy conjecture
predicts that, if $\alpha\in \Q$ is a pole of
$S(\mathscr{F}_x,\phi/dt;\LL^{-s})$, then $\exp(2\pi i \alpha)$ is
a monodromy eigenvalue of $f$ at some closed point of the special
fiber of $f$. The conjecture has been solved, for instance, in the
case where $X$ is a surface \cite{Loepadic}\cite{Rod}, but the
general case remains wide open.

We will formulate a global form of this conjecture, and prove it
for abelian varieties. Denote by $\tau(c(A))$ the order of $c(A)$
in the group $\Q/\Z$, and by $\Phi_{\tau(c(A))}(t)$ the cyclotomic
polynomial whose roots are the primitive roots of unity of order
$\tau(c(A))$. Let $\sigma$ be a topological generator of the tame
monodromy group $G(K^t/K)$. We give a formula for the
characteristic polynomial of the action of $\sigma$ on the Tate
module of $A$, in terms of the jumps in Edixhoven's filtration. We
also prove that $\Phi_{\tau(c(A))}(t)$ divides the characteristic
polynomial of $\sigma$ on $H^g(A\times_K K^t,\Q_\ell)$, with $g$
the dimension of $A$.  Since $s=c(A)$ is the only pole of
$Z_A(\LL^{-s})$, this yields a global form of Denef and Loeser's
monodromy conjecture for abelian varieties.

Of course, we hope that our proofs and results will lead to new
insights into the local case of the conjecture. In this respect,
it is intriguing that the Greenberg schemes appearing in the
construction of Edixhoven's filtration also play a fundamental
role in the theory of motivic integration and local zeta functions
\cite{NiSe}.

Although the main results in this paper concern abelian varieties,
we develop a large part of the theory in greater generality, in
order to include the case of semi-abelian varieties. It would be
 interesting to know if the monodromy conjecture holds also for tamely ramified semi-abelian varieties.
 The case of
algebraic tori is treated in the follow-up paper \cite{Ni-tori}.
The assumption that $A$ is tamely ramified is crucial for our
arguments; it is a challenging problem to adapt our results to
wildly ramified (semi-)abelian varieties.


\bigskip
To conclude this introduction, we give a survey of the structure
of the paper. Section \ref{sec-prelim} gathers some preliminaries
on motivic integration and motivic generating series. In
Definition \ref{def-gmp}, we formulate a global version of Denef
and Loeser's monodromy conjecture.

Section \ref{sec-neron} deals with N\'eron models of smooth
commutative algebraic $K$-groups $G$. The N\'eron models we
consider are the maximal quasi-compact open subgroups of the
locally of finite type N\'eron model $\mathcal{G}^{lft}$ of $G$.
 The proof of their existence boils down to
showing that the component group of $\mathcal{G}^{lft}_s$ is
finitely generated (Proposition \ref{lemma-finite}).
 We
characterize the N\'eron model by a universal property in
Definition \ref{def-neron}.

Section \ref{sec-edixchai} contains the basic results on jumps,
the base change conductor, and the relation between them. In
Section \ref{subsec-edix} we extend Edixhoven's results in
\cite{edix} to arbitrary smooth commutative algebraic $K$-groups
that admit a N\'eron model, and in particular to semi-abelian
varieties. Section \ref{subsec-chai} briefly reviews Chai's
elementary divisors and base change conductor for semi-abelian
$K$-varieties, and in Section \ref{subsec-compar} we show that for
tamely ramified semi-abelian $K$-varieties, the jumps and the
elementary divisors are equivalent (Corollary \ref{cor-compar}).

In Section \ref{sec-monodromy} we study the relation between jumps
and monodromy eigenvalues for tamely ramified abelian
$K$-varieties $A$. Theorem \ref{theo-monodromy} computes the
characteristic polynomial of the tame monodromy operator $\sigma$
on the Tate module $T_\ell A$ in terms of the jumps of $A$. In
Corollary \ref{cor-monocond}, we prove that $c(A)$ corresponds to
a monodromy eigenvalue on $H^g(A\times_K K^t,\Q_\ell)$ as
explained above.

Section \ref{sec-tamebc} deals with the behaviour of
$\mathcal{A}(d)^o_s$ when $d$ varies over $\N'$. In Section
\ref{sec-ord} we express $ord_{\mathcal{A}(d)_s^o}(\omega(d))$ in
terms of the jumps of $A$. In Section \ref{sec-rat}, we define the
motivic zeta function of a semi-abelian $K$-variety, and we prove
the global monodromy conjecture for tamely ramified abelian
varieties (Theorem \ref{thm-ratzeta}).

\section{Preliminaries}\label{sec-prelim}
\subsection{Notation}
We denote by $R$ a discrete valuation ring, by $K$ its quotient
field, and by $k$ its residue field.  Additional conditions on $R$
and $k$ will be indicated at the beginning of each section. If $R$
has equal characteristic, then we fix a $k$-algebra structure on
$R$ such that the composition $k\rightarrow R\rightarrow k$ is the
identity. We denote by $p$ the characteristic exponent of $k$, and
we fix a uniformizer $\pi$ in $R$. The choice of a value $|\pi|$
in $]0,1[$ determines a $\pi$-adic absolute value $|\cdot|$ on
$K$. We denote by $\N'$ the set of strictly positive integers that
are prime to $p$, and we fix a prime $\ell\neq p$. We fix a
separable closure $K^s$ of $K$, and we denote by $K^t$ the tame
closure of $K$ in $K^s$.

We recall the following definition \cite[3.6.1]{neron}. A morphism
of discrete valuation rings $R\rightarrow R'$ is called unramified
if it is a flat local morphism, $\pi$ is a uniformizer in $R'$,
and the extension of residue fields $R/(\pi)\rightarrow R'/(\pi)$
is separable. In this case, we call $R'$ an unramified extension
of $R$.

We denote by $K_0(Var_k)$ the Grothendieck ring of $k$-varieties,
by $\LL$ the class $[\A^1_k]$ of the affine line $\A^1_k$ in
$K_0(Var_k)$, and by $\mathcal{M}_k$ the localization of
$K_0(Var_k)$ w.r.t. $\LL$. See for instance
\cite[2.1]{Ni-tracevar}. We denote by
$$\chi_{top}:\mathcal{M}_k\rightarrow \Z$$ the unique ring
morphism that sends the class $[X]$ of every $k$-variety $X$ to
the $\ell$-adic Euler characteristic $\chi_{top}(X)$ of $X$. This
morphism is independent of $\ell$.


For every scheme $S$, we denote by $(Sch/S)$ the category of
$S$-schemes. If $A$ is a commutative ring, then we write $(Sch/A)$
instead of $(Sch/\Spec A)$. We denote by
$$(\cdot)_K:(Sch/R)\rightarrow (Sch/K):X\mapsto X_K=X\times_R K$$ the generic
fiber functor, and by
$$(\cdot)_s:(Sch/R)\rightarrow (Sch/k):X\mapsto X_s=X\times_R k$$ the special
fiber functor. For every $R$-scheme $X$ and every section $\psi$
in $X(R)$, we denote by $\psi(0)$ the image in $X_s$ of the closed
point of $\Spec R$.

 For every scheme $S$ and every  group $S$-scheme $G$, we denote by
$e_G\in G(S)$ the unit section, and we put
$\omega_{G/S}=e_{G}^*\Omega^1_{G/S}$. If $F$ is a field, then an
algebraic $F$-group is a group $F$-scheme of finite type. We call
a semi-abelian $K$-variety tamely ramified if it acquires
semi-abelian reduction on a finite tame extension of $K$.

For every real number $x$ we denote by $\lfloor x\rfloor$ the
unique integer in $]x-1,x]$ and by $\lceil x\rceil$ the unique
integer in $[x,x+1[$. We put $[x]=x-\lfloor x \rfloor \in [0,1[$.
We denote by
$$\tau:\Q\rightarrow \Z_{>0}$$ the map that sends a rational
number $a$ to its order $\tau(a)$ in the quotient group $\Q/\Z$.
For each $i\in \Z_{>0}$, we denote by $\Phi_i(t)\in \Z[t]$ the
cyclotomic polynomial whose zeroes are the primitive $i$-th roots
of unity.

If $K$ is strictly henselian, then we adopt the following
notations.
 We fix a topological generator $\sigma$
of the tame monodromy group $G(K^t/K)$. For each $d\in \N'$ we
denote by $K(d)$ the unique degree $d$ extension of $K$ in $K^t$,
and by $R(d)$ the normalization of $R$ in $K(d)$. For every
algebraic $K$-variety $X$ we put $X(d)=X\times_K K(d)$. For every
differential form $\omega$ on $X$ we denote by $\omega(d)$ the
pull-back of $\omega$ to $X(d)$.

Finally, assume that $R$ is complete. We call a formal $R$-scheme
$\mathfrak{Y}$ $stft$ if it is separated and topologically of
finite type over $R$. We denote by $\mathfrak{Y}_s$ its special
fiber (a separated $k$-scheme of finite type) and by
$\mathfrak{Y}_\eta$ its generic fiber (a separated quasi-compact
rigid $K$-variety). If $k$ is separably closed, $X$ is a rigid
$K$-variety and $\omega$ a differential form on $X$, then $X(d)$
and $\omega(d)$ are defined as in the algebraic case. We denote by
$(\cdot)^{an}$ the rigid analytic GAGA functor.


\subsection{Order of a gauge form}
Let $X$ be a smooth $R$-scheme of pure relative dimension $g$, and
let $\omega$ be a gauge form on $X_K$, i.e. a nowhere vanishing
differential $g$-form on $X_K$.  Let $C$ be a connected component
of $X_s$, and denote by $\eta_C$ the generic point of $C$. Then
$\mathcal{O}_{X,\eta_C}$ is a discrete valuation ring, whose
maximal ideal is generated by $\pi$. Let $\psi$ be a section in
$X(R)$. Recall the following definitions.
\begin{definition}[Order of a gauge form]\label{def-ord}
Choose an element $\alpha$ in $\N$ such that  $\pi^{\alpha}\omega$
 belongs to the image of the natural injection
$$\Omega^g_{X/R}(X)\rightarrow \Omega^g_{X_K/K}(X_K)$$
 We denote by $\omega'$ the unique inverse image
of $\pi^{\alpha}\omega$ in $\Omega^g_{X/R}(X)$.

 The order of $\omega$ along $C$ is defined by
 $$ord_{C}(\omega)=\mathrm{length}_{\mathcal{O}_{X,\eta_C}}\left((\Omega^g_{X/R})_{\eta_C}/\mathcal{O}_{X,\eta_C}\cdot \omega'
 \right)-\alpha$$
 The order of $\omega$ at $\psi$ is defined by
 $$ord(\omega)(\psi)=\mathrm{length}_{R}\psi^*\left(\Omega^g_{X/R}/\mathcal{O}_{X}\cdot
 \omega'\right)-\alpha$$
 These definitions are independent of the choices of $\pi$ and
 $\alpha$.
\end{definition}
Note that $ord_{C}(\omega)$ and $ord(\omega)(\psi)$ are finite,
since $\omega$ is a gauge form.

\begin{prop}\label{prop-ordsec}
Let $X$ be a smooth $R$-scheme of pure relative dimension, $C$ a
connected component of $X_s$, and $\psi$ a section in $X(R)$ such
that $\psi(0)\in C$. For every gauge form $\omega$ on $X_K$, we
have
$$ord_{C}(\omega)=ord(\omega)(\psi)$$
\end{prop}
\begin{proof}
This property follows immediately from its formal counterpart
\cite[5.10]{ni-trace} and can easily be proven directly by similar
arguments.
\end{proof}

Now assume that $R$ is complete. If $\X$ is a smooth $stft$ formal
$R$-scheme, then the order of a gauge form $\phi$ on $\X_\eta$
along a connected component $C$ of the special fiber $\X_s$ is
defined as in the algebraic case \cite[4.3]{formner}. If $\X$ is
the formal $\pi$-adic completion of $X$, and $\phi$ is the
restriction to $\X_\eta$ of the rigid analytification
$\omega^{an}$ on $(X_K)^{an}$, then we have
$$ord_C(\omega)=ord_C(\phi)$$

\subsection{Motivic integration on rigid varieties}
Assume that $K$ is complete and that $k$ is perfect. Let $X$ be a
separated smooth quasi-compact rigid $K$-variety of pure
dimension, and let $\omega$ be a gauge form on $X$. The motivic
integral
$$\int_{X}|\omega|\quad \in \mathcal{M}_k$$ was defined in
\cite[4.1.2]{motrigid}. For our purposes, the following
proposition can serve as a definition. The result is merely
slightly stronger than \cite[4.3.1]{motrigid}, but the difference
is important for the applications in this paper. Recall that a
weak N\'eron model for $X$ is a smooth $stft$ formal $R$-scheme
$\mU$ endowed with an open immersion of rigid $K$-varieties
$i:\mU_\eta\rightarrow X$ such that $i(K'):\mU_\eta(K')\rightarrow
X(K')$ is bijective for every finite unramified extension $K'$ of
$K$ \cite[1.3]{formner}. Such a weak N\'eron model always exists,
by \cite[3.3]{formner}.
\begin{prop}\label{prop-motint}
For every weak N\'eron model $(\mU,i)$ of $X$, we have
\begin{equation}\label{eq-motint}
\int_{X}|\omega|=\LL^{-dim(X)}\sum_{C\in
\pi_0(\mU_s)}[C]\LL^{-ord_C (i^*\omega)}\quad \in
\mathcal{M}_k\end{equation}
 where $\pi_0(\mU_s)$ denotes the set
of connected components of $\mU_s$. In particular, the right hand
side of (\ref{eq-motint}) does not depend on $\mU$.
\end{prop}
\begin{proof}
This is almost the statement of \cite[4.3.1]{motrigid}, except
that there the additional condition was imposed that $\mU$ is an
open formal subscheme of a formal model of $X$. Let us explain why
it can be omitted. By \cite[4.3.1]{motrigid}, the right hand side
of (\ref{eq-motint}) equals $\int_{\mU_\eta}|i^*\omega|$. Hence,
it suffices to show that
$$\int_{\mU_\eta}|i^*\omega|=\int_{X}|\omega|$$
This follows from \cite[5.9]{NiSe-weilres}.
\end{proof}
If $X$ does not have pure dimension, then by a gauge form $\omega$
on $X$, we mean the datum of a gauge form $\omega_Y$ on each
connected component $Y$ of $X$. The motivic integral of $\omega$
on $X$ is then defined by
$$\int_{X}|\omega|=\sum_{Y\in \pi_0(X)}\int_{Y}|\omega_Y|\quad \in \mathcal{M}_k$$
\subsection{Motivic generating series}\label{subsec-motseries}
Assume that $R$ is complete and that $k$ has characteristic zero.
Let $X$ be a separated smooth quasi-compact rigid $K$-variety, and
$\omega$ a gauge form on $X$. We denote by $S(X,\omega;T)$ the
motivic generating series
$$S(X,\omega;T)
\in \mathcal{M}_k[[T]]$$ from \cite[7.2]{NiSe}. It depends on the
choice of uniformizer $\pi$ in general, but it is independent of
this choice if $k$ is algebraically closed \cite[4.10]{ni-trace}.
In fact, if $k$ is algebraically closed, then
$$S(X,\omega;T)=\sum_{d>0}\left(\int_{X(d)}|\omega(d)|\right)T^d\quad
\in \mathcal{M}_k[[T]]$$

In any case, the series $S(X,\omega;T)$ can be computed explicitly
on a regular $stft$ formal $R$-model $\X$ of $X$ whose special
fiber $\X_s$ is a divisor with strict normal crossings
\cite[7.7]{NiSe}. Such a model always exists, by embedded
resolution of singularities for generically smooth $stft$ formal
$R$-schemes \cite[3.4.1]{temkin-resol}. The explicit expression
for the motivic generating series $S(X,\omega;T)$ shows in
particular that $S(X,\omega;T)$ is rational and belongs to the
subring
$$\mathcal{M}_k\left[\frac{\LL^a T^b}{1-\LL^a
T^b}\right]_{(a,b)\in \Z\times \Z_{>0}}$$ of $\mathcal{M}_k[[T]]$.

\subsection{The monodromy conjecture} In this section we assume
that $k$ has characteristic zero and that $R$ is complete.

\subsubsection*{The local case.}
If $\mathfrak{X}$ is a regular $stft$ formal $R$-scheme, then its
motivic Weil generating series
$$S(\mathfrak{X};T)\in \mathcal{M}_{\X_s}[[T]]$$ was defined in
\cite[7.33]{ni-trace}. Here $\mathcal{M}_{\X_s}$ is the localized
Grothendieck ring of $\mathcal{M}_{\X_s}$-varieties; see
\cite[2.1]{Ni-tracevar}. The image of $S(\X;T)$ under the
forgetful morphism
$$\mathcal{M}_{\X_s}[[T]]\rightarrow \mathcal{M}_k[[T]]$$ equals
$S(\X_\eta,\omega;T)$, with $\omega$ a so-called
\textit{Gelfand-Leray form} on $\X_\eta$ \cite[7.21]{ni-trace}.
\begin{conjecture}[Local Motivic Monodromy
Conjecture]\label{conj-lmc} If $\X$ is a regular $stft$ formal
$R$-scheme, then there exists a finite subset $\mathcal{S}$ of
$\Z\times \Z_{>0}$ such that
$$S(\X;T)\in \mathcal{M}_{\X_s}\left[T,\frac{1}{1-\LL^a T^b}\right]_{(a,b)\in \mathcal{S}}$$ and such that
for each $(a,b)\in \mathcal{S}$, the cyclotomic polynomial
$\Phi_{\tau(a/b)}(t)$ divides the characteristic polynomial of
$\sigma$ on $R^i\psi_{\X}(\Q_\ell)_x$, for some $i\in \Z_{\geq 0}$
and some geometric closed point $x$ of $\X_s$.
\end{conjecture}

 Here $\psi_{\X}$ denotes the nearby cycle functor on the formal
 $R$-scheme $\X$ (called vanishing cycle functor in \cite{Berk-vanish}\cite{berk-vanish2}). In particular, Conjecture
\ref{conj-lmc} implies that for each pole $\alpha$ of
$S(\mathfrak{X};\LL^{-s})$ and for every embedding $\Q_\ell
\hookrightarrow \C$, the value $\exp(2\pi i \alpha)$ is an
eigenvalue of the monodromy action of $\sigma$ on
$R^i\psi_{\X}(\Q_\ell)_x$, for some $i\in \Z_{\geq 0}$ and some
geometric closed point $x$ of $\X_s$. By the comparison result in
\cite[9.6]{ni-trace}, Conjecture \ref{conj-lmc} implies the
following conjecture of Denef and Loeser's, which at its turn
generalizes Igusa's monodromy conjecture for $p$-adic zeta
functions \cite[\S\,2.4]{DL5}.

\begin{conjecture}[Denef and Loeser's Motivic Monodromy
Conjecture]\label{conj-dl} Let $k$ be a subfield of $\C$, and let
$X$ be a smooth and irreducible $k$-variety endowed with a
dominant morphism $f:X\rightarrow \A^1_k$. Denote by $X_s$ the
special fiber of $f$ and by
$$Z_f(T)\in \mathcal{M}_{X_s}[[T]]$$
the motivic zeta function associated to $f$ \cite[3.2.1]{DL3}.
Then there exists a finite subset $\mathcal{S}$ of $\Z\times
\Z_{>0}$ such that
$$Z_f(T)\in \mathcal{M}_{X_s}\left[T,\frac{1}{1-\LL^a T^b}\right]_{(a,b)\in \mathcal{S}}$$ and such that
for each $(a,b)\in \mathcal{S}$, the value $\exp (2\pi i a/b)$ is
an eigenvalue of monodromy on $R^i\psi_f(\Q)_x$ for some $i\in
\Z_{\geq 0}$ and some $x\in X_s(\C)$. Here $R\psi_f(\Q)\in
D^b_c(X_s(\C),\Q)$ denotes the complex-analytic nearby cycles
complex of $f$.
\end{conjecture}
\noindent We refer to \cite{DenefBour}\cite{Ni-japan}\cite{NiSe3}
for an introduction to the $p$-adic and motivic monodromy
conjecture.

\begin{remark}
To be precise, Denef and Loeser's conjecture is a bit stronger
than Conjecture \ref{conj-dl}, because it is stated for the
``monodromic'' motivic zeta function which carries an additional
action of the pro-$k$-group $\widehat{\mu}$ of roots of unity
\cite[3.2.1]{DL3}. The zeta function $Z_f(T)$ in Conjecture
\ref{conj-dl} is the image of this monodromic zeta function under
the forgetful morphism
$$\mathcal{M}^{\widehat{\mu}}_{X_s}[[T]]\rightarrow
\mathcal{M}_{X_s}[[T]]$$
\end{remark}

\subsubsection*{The global case}
 The following definition formulates a global
version of the motivic monodromy conjecture.
\begin{definition}\label{def-gmp}
Let $X$ be a smooth, proper, geometrically connected $K$-variety,
and assume that $X$ admits a gauge form $\omega$. We say that $X$
satisifes the Global Monodromy Property (GMP) if there exists a
finite subset $\mathcal{S}$ of $\Z\times \Z_{>0}$ such that
$$S(X,\omega,T)\in \mathcal{M}_k\left[T,\frac{1}{1-\LL^a T^b}\right]_{(a,b)\in \mathcal{S}}$$ and such that for each
$(a,b)\in \mathcal{S}$, the cyclotomic polynomial
$\Phi_{\tau(a/b)}(t)$ divides the characteristic polynomial of
$\sigma$ on $H^i(X\times_K K^s,\Q_\ell)$ for some $i\in \N$.
\end{definition}

\noindent Note that the property $(GMP)$ only depends on $X$, and
not on $\omega$, since we have
$$S(X,u\cdot \omega;T)=S(X,\omega;\LL^{-v_K(u)}T)\in
\mathcal{M}_k[[T]]$$ for all $u \in K^*$, where $v_K$ denotes the
discrete valuation on $K^*$.

\begin{question}
Is there a natural condition on $X$ that guarantees that $X$
satisfies the Global Monodromy Property (GMP)?
\end{question}

\noindent We will show in Theorem \ref{thm-ratzeta} that if $k$ is
algebraically closed, every abelian $K$-variety satisfies the
Global Monodromy Property. Moreover, we will show that this result
extends to tamely ramified abelian varieties in mixed and positive
characteristic. Our proof does not use resolution of
singularities.

%
%
%
%
%

\section{N\'eron models}\label{sec-neron}
Let $G$ be a smooth commutative algebraic $K$-group. The notion of
 N\'eron $lft$-model for $G$ over $R$ is defined in
\cite[10.1.1]{neron}. It is shown in \cite[10.2.2]{neron} that $G$
admits a N\'eron $lft$-model iff $G\times_K \widehat{K}$ does not
contain a subgroup of type $\mathbb{G}_{a,\widehat{K}}$, where we
denote by $\widehat{K}$ the completion of $K$. In particular,
every semi-abelian $K$-variety $A$ admits a N\'eron $lft$-model
over $R$.


\begin{lemma}\label{lemma-conn}
Assume that $k$ is separably closed. Let $G$ be a smooth
commutative algebraic $K$-group. Then $G$ admits a N\'eron
$lft$-model $\mathcal{G}^{lft}$ iff $G^o$ admits a N\'eron
$lft$-model $\mathcal{H}^{lft}$. If we denote by $f$ the unique
morphism $\mathcal{H}^{lft}\rightarrow \mathcal{G}^{lft}$ of group
$R$-schemes extending the open immersion $G^o\rightarrow G$, then
the morphism of constant groups
$\pi_0(f_s):\pi_0(\mathcal{H}^{lft}_s)\rightarrow
\pi_0(\mathcal{G}^{lft}_s)$ is an injection whose cokernel is
canonically isomorphic to $G(K)/G^o(K)$.
\end{lemma}
\begin{proof}
It is clear from the definition that $G$ admits a N\'eron
$lft$-model iff $G^o$
 does, and that the formation of N\'eron
$lft$-models commutes with finite disjoint unions, so that $f$ is
an open and closed immersion. If $C$ is a connected component of
$G$ without $K$-point, then $C$ is a N\'eron $lft$-model of
itself, with empty special fiber.
 The connected
components of $G$ with $K$-point form a constant subgroup scheme
of $\pi_0(G)$, which is canonically isomorphic to $G(K)/G^o(K)$.
For any connected component $C$ of $\mathcal{G}^{lft}_s$, there
exists a unique connected component $C'$ of $G$ such that $C$
belongs to the schematic closure of $C'$ in $\mathcal{G}^{lft}$,
and $C'(K)\neq \emptyset$. The map $C\mapsto C'$ defines a
surjection $\pi_0(\mathcal{G}^{lft}_s)\rightarrow G(K)/G^o(K)$
whose kernel is precisely $\pi_0(\mathcal{H}^{lft}_s)$.
\end{proof}
\begin{lemma}\label{lemma-sep}
Assume that $K$ is Henselian. Let $G$ be a smooth commutative
algebraic $K$-group, and let $K'$ be a finite separable extension
of $K$. Then $G$ admits a N\'eron $lft$-model iff $G\times_K K'$
admits a N\'eron $lft$-model.
\end{lemma}
\begin{proof}
By Lemma \ref{lemma-conn} we may assume that $G$ is connected. By
\cite[10.2.2]{neron} it suffices to show that $G\times_K
\widehat{K}$ admits a subgroup of type
$\mathbb{G}_{a,\widehat{K}}$ iff $G\times_K \widehat{K'}$ admits a
subgroup of type $\mathbb{G}_{a,\widehat{K'}}$. The ``only if''
part is obvious, so let us prove the ``if'' part. Since
$\widehat{K'}$ is a finite separable extension of $\widehat{K}$,
we may assume that $K$ is complete. Assume that $G$ does not admit
a subgroup of type $\mathbb{G}_{a,K}$. Let $T$ be the maximal
torus in $G$, and put $H=G/T$. It is shown in the proof of
\cite[10.2.2]{neron} that $H\times_K K'$ does not admit a subgroup
of type $\mathbb{G}_{a,K'}$. Assume that $G\times_K K'$ has a
subgroup isomorphic to $\mathbb{G}_{a,K'}$. Then the restriction
of $f:G\times_K K'\rightarrow H\times_K K'$ to $\mathbb{G}_{a,K'}$
must have a non-trivial kernel $U$, which is contained in
$ker(f)=T\times_K K'$. But $U$ is unipotent, and there are no
non-trivial morphisms of group $K'$-schemes $U\rightarrow
T\times_K K'$ \cite[XVII.2.4]{sga3.2}, so we arrive at a
contradiction.
%
\end{proof}
\begin{lemma}\label{lemma-ft}
Assume that $K$ is Henselian. Let $K'$ be a finite separable
extension of $K$, and denote by $R'$ the normalization of $R$ in
$K'$. Let $G$ be a smooth commutative algebraic $K$-group that
admits a N\'eron $lft$ model $\mathcal{G}^{lft}$, and denote by
$(\mathcal{G}')^{lft}$ the N\'eron $lft$-model of $G'=G\times_K
K'$. Consider the unique morphism of group $R'$-schemes
$$f:\mathcal{G}^{lft}\times_R R'\rightarrow (\mathcal{G}')^{lft}$$
that extends the canonical isomorphism between the generic fibers.
This morphism $f$ is of finite type.
\end{lemma}
\begin{proof}
 Consider the morphism of group $R$-schemes
$$g:\mathcal{G}^{lft}\rightarrow
W:=\prod_{R'/R}(\mathcal{G}')^{lft}$$ obtained from $f$ by
adjunction, and denote by $\mathcal{H}$ the schematic closure of
$G_K$ in $W_K$. It is a closed subgroup scheme of $W$. Since
$\mathcal{G}^{lft}$ is flat, the morphism $g$ factors through a
morphism of group $R$-schemes $g':\mathcal{G}^{lft}\rightarrow
\mathcal{H}$. By \cite[10.1.4]{neron}, this is the canonical group
smoothening of $\mathcal{H}$. In particular, it is of finite type.
This implies that $f$ is of finite type, since it is the
composition of $g\times_R R'$ with the tautological morphism of
finite type $W\times_R R'\rightarrow (\mathcal{G}')^{lft}$.
\end{proof}
\begin{lemma}\label{lemma-splittor}
Let $G$ be a connected smooth commutative algebraic $K$-group that
admits a N\'eron $lft$ model $\mathcal{G}^{lft}$. Assume that the
maximal torus $T$ in $G$ is split. Then the component group
$\pi_0(\mathcal{G}^{lft}_s\times_k k^s)$ is finitely generated.
\end{lemma}
\begin{proof}
Since the formation of N\'eron $lft$-models commutes with
unramified base change \cite[10.1.3]{neron} we may assume that $K$
is complete and strictly henselian. Denote the N\'eron $lft$-model
of $T$ by $\mathcal{T}^{lft}$, and denote by $H$ the quotient
$G/T$. Then $H$ has a N\'eron $lft$-model $\mathcal{H}^{lft}$
which is of finite type over $R$, by the proof of
\cite[10.2.2]{neron}.
 By the proof of \cite[10.1.7]{neron}, the sequence
$$\pi_0(\mathcal{T}_s^{lft})\rightarrow \pi_0(\mathcal{G}^{lft}_s)\rightarrow \pi_0(\mathcal{H}^{lft}_s)\rightarrow 1$$
is exact. Since $\pi_0(\mathcal{T}^{lft}_s)$ is a free $\Z$-module
whose rank is equal to the dimension of $T$ and
$\pi_0(\mathcal{H}^{lft}_s)$ is finite, we see that
$\pi_0(\mathcal{G}^{lft}_s)$ is finitely generated.
\end{proof}
\begin{prop}\label{lemma-finite}
 Let $G$ be a smooth
commutative algebraic $K$-group, and assume that $G$ admits a
N\'eron $lft$-model $\mathcal{G}^{lft}$. Then the component group
$\pi_0(\mathcal{G}^{lft}_s\times_k k^s)$ is finitely generated. In
particular, its torsion part is finite.
\end{prop}
\begin{proof}
%
Since the formation of N\'eron $lft$-models commutes with
unramified base change \cite[10.1.3]{neron} we may assume that $K$
is strictly Henselian. We may also assume that $G$ is connected,
by Lemma \ref{lemma-conn}. Let $T$ be the maximal torus in $G$.
Let $K'$ be a finite separable extension of $K$ such that $T$
splits over $K'$. Denote by $R'$ the normalization of $R$ in $K'$.
By Lemma \ref{lemma-sep}, $G\times_K K'$ admits a N\'eron
$lft$-model $(\mathcal{G}')^{lft}$. By Lemma \ref{lemma-splittor},
the component group $\pi_0((\mathcal{G}')^{lft}_s)$ is finitely
generated. By Lemma \ref{lemma-ft}, the natural morphism of group
$R'$-schemes
$$f:\mathcal{G}^{lft}\times_R R'\rightarrow (\mathcal{G}')^{lft}$$
is of finite type, so that the induced morphism of component
groups
$$\pi_0(\mathcal{G}^{lft}_s)\rightarrow
\pi_0((\mathcal{G}')^{lft}_s)$$ has finite kernel.
Hence, $\pi_0(\mathcal{G}^{lft}_s)$ is finitely generated. 
\end{proof}

\begin{definition}\label{def-neron}
For every commutative ring $S$, we denote by $(\mathscr{G}r/S)$
the category of smooth separated group $S$-schemes of finite type.
We denote by $(Groups)$ the category of groups.

 Let $G$ be a smooth algebraic $K$-group. If the
functor
$$NM_G:(\mathscr{G}r/R)\rightarrow (Groups):H\mapsto
Hom_{(\mathscr{G}r/K)}(H_K,G)$$ is representable by an object
$\mathcal{G}$ of $(\mathscr{G}r/R)$, then we call $\mathcal{G}$
the N\'eron model of $G$. We denote by $\Phi_G$ the group
$k$-scheme of connected components $\pi_0(\mathcal{G}_s)$, and we
denote its rank by $\phi(G)$.
\end{definition}
It follows immediately from the definition that there exists a
canonical isomorphism of algebraic $K$-groups $G\cong
\mathcal{G}_K$. If $G$ admits a N\'eron model in the sense of
\cite{neron} (i.e., if its N\'eron $lft$-model exists and is
quasi-compact) then this N\'eron model represents the functor
$NM_G$, so that our definition includes the one in \cite{neron}.

\begin{prop}\label{prop-lftner}
Let $G$ be a smooth commutative algebraic $K$-group, and assume
that $G$ has a N\'eron $lft$-model $\mathcal{G}^{lft}$. Then $G$
has a N\'eron model $\mathcal{G}$.
 The canonical isomorphism of group $K$-schemes
$\mathcal{G}_K\cong G$ extends uniquely to a morphism of group
$R$-schemes
$$\varphi:\mathcal{G}\rightarrow \mathcal{G}^{lft}$$
 The morphism $\varphi$ is an
open immersion, and its image is the maximal quasi-compact open
subgroup scheme of $\mathcal{G}^{lft}$. The group $k$-scheme
$\Phi_G$ is the torsion part of $\pi_0(\mathcal{G}^{lft}_s)$.
\end{prop}
\begin{proof}
Denote by $\mathcal{G}$ the union of the generic fiber $G$ of
$\mathcal{G}^{lft}$ with all connected components $C$ of the
special fiber $\mathcal{G}^{lft}_s$ such that $C$ defines a
torsion point of
 $\pi_0(\mathcal{G}^{lft}_s)$. Then
$\mathcal{G}$ is an open subgroup scheme of $\mathcal{G}^{lft}$.
Moreover, since $\pi_0(\mathcal{G}^{lft}_s\times_k k^s)$ has
finite torsion part by Proposition \ref{lemma-finite},
$\mathcal{G}$ is quasi-compact and hence of finite type over $R$.
It is clear that $\mathcal{G}$ is the largest quasi-compact open
subgroup scheme of $\mathcal{G}^{lft}$.

Now we show that $\mathcal{G}$ represents the functor $NM_G$. Let
$H$ be a smooth group scheme of finite type over $R$, and consider
a morphism of group $K$-schemes $H_K\rightarrow G$. By the
universal property of the N\'eron $lft$-model, this morphism
extends uniquely to a morphism of group $R$-schemes
$h:H\rightarrow \mathcal{G}^{lft}$. Since $H$ is quasi-compact,
the image of $\pi_0(H_s\times_k k^s)$ in
$\pi_0(\mathcal{G}_s^{lft}\times_k k^s)$ is a finite subgroup of
$\pi_0(\mathcal{G}_s^{lft}\times_k k^s)$. Hence, the image of $h$
is contained in $\mathcal{G}$.
\end{proof}

\begin{cor}
If $A$ is a semi-abelian $K$-variety, then $A$ admits a N\'eron
model.
\end{cor}

\begin{definition}\label{def-bounded}
Assume that $K$ is complete. Let $G$ be a smooth commutative
algebraic $K$-group that admits a N\'eron model $\mathcal{G}$. We
define the bounded part $G^b$ of $G$ as the generic fiber of the
formal $\pi$-adic completion of $\mathcal{G}$. It is a separated,
smooth, quasi-compact open rigid subgroup of the rigid
analytification $G^{an}$ of $G$.
\end{definition}


\begin{prop}\label{prop-equiv}
Assume that $R$ is excellent. Let $G$ be a smooth commutative
algebraic $K$-group. Then $G$ admits a N\'eron model iff $G$
admits a N\'eron $lft$-model.
\end{prop}
\begin{proof}
If $G$ admits a N\'eron $lft$-model, then $G$ admits a N\'eron
model by Proposition \ref{prop-lftner}. Conversely, assume that
$G$ admits a N\'eron model $\mathcal{G}$. By \cite[10.2.2]{neron}
it is enough to show that $G$ does not contain a subgroup of type
$\mathbb{G}_{a,K}=\Spec K[\xi]$.

Suppose that $G$ admits a subgroup of type $\mathbb{G}_{a,K}$.
Since any subgroup scheme of $G$ is closed
\cite[VI$_B$\,1.4.2]{sga3.1}, the set
$$B=\mathcal{G}(R)\cap \mathbb{G}_{a,K}(K)\subset \mathbb{G}_{a,K}(K)$$
is bounded in $\mathbb{G}_{a,K}$ in the sense of
\cite[1.1.2]{neron}. Therefore, $B$ is contained in a closed disc
in $\mathbb{G}_{a,K}(K)$ defined by $|\xi|\leq |\pi|^{-N}$ for
some $N\in \N$.

Now consider the smooth group $R$-scheme of finite type
$$H=\mathbb{G}_{a,R}=\Spec R[\zeta]$$
and the morphism of group $K$-schemes
$$h:H_K\rightarrow \mathbb{G}_{a,K}:\xi\mapsto \pi^{-(N+1)}\zeta$$
By the universal property of the N\'eron model, the induced
morphism $H_K\rightarrow G$ extends to an $R$-morphism
$H\rightarrow \mathcal{G}$. This is a contradiction, since the
point $x$ of $H_K(K)$ defined by the ideal $(\zeta-1)$ belongs to
$H(R)$, while it is mapped to the point $(\xi-\pi^{-(N+1)})$ in
$\mathbb{G}_{a,K}(K)$, which does not belong to $B$.
\end{proof}


%
%

\begin{prop}\label{prop-unram}
Let $R\rightarrow S$ be an unramified morphism of discrete
valuation rings and denote by $L$ the quotient field of $S$. Let
$G$ be a smooth commutative algebraic $K$-group. If $G$ admits a
N\'eron $lft$-model $\mathcal{G}^{lft}$, and $\mathcal{G}$ is a
N\'eron model of $G$, then $\mathcal{G}\times_R S$ is a N\'eron
model for $G\times_K L$.
\end{prop}
\begin{proof}
By \cite[10.1.3]{neron},  $\mathcal{G}^{lft}\times_R S$ is a
N\'eron $lft$-model for $G\times_K L$. Since the torsion part of
$\pi_0(\mathcal{G}^{lft}_s)$ is stable under the base change
$k\rightarrow S/(\pi)$, it follows from Proposition
\ref{prop-lftner} that $\mathcal{G}\times_R S$ is a N\'eron model
of $G\times_K L$.
\end{proof}

\begin{prop}\label{prop-maxbound}
Assume that $R$ is excellent. Let $G$ be a smooth commutative
algebraic $K$-group, and assume that $G$ has a N\'eron model
$\mathcal{G}$. Let $R\rightarrow S$ be an unramified morphism of
discrete valuation rings and denote by $L$ the quotient field of
$S$. Then $\mathcal{G}(S)\subset G(L)$
 is the unique maximal subgroup of $G(L)$
that is bounded in $G$ in the sense of \cite[1.1.2]{neron}.
\end{prop}
\begin{proof}
By Propositions \ref{prop-equiv} and \ref{prop-unram}, and
\cite[1.1.5]{neron}, we may assume that $R=S$. The group
$\mathcal{G}(R)$ is bounded in $G$ by \cite[1.1.7]{neron}.
Conversely, let $x$ be a point of $G(K)$ and assume that $x$
belongs to a subgroup of $G(K)$ that is bounded in $G$. Then the
subgroup $<\!\!x\!\!>$ of $G(K)$ generated by $x$ is bounded in
$G$. We will show that $x\in \mathcal{G}(R)$. By Proposition
\ref{prop-equiv}, $G$ has a N\'eron $lft$-model
$\mathcal{G}^{lft}$. The point $x$ extends uniquely to a section
$\psi$ in $\mathcal{G}^{lft}(R)$, and it suffices to show that
$\psi_s$ is contained in $\mathcal{G}_s(k)\subset
\mathcal{G}^{lft}_s(k)$.

Since $<\!\!x\!\!>$ is bounded in $G$, there exists a smooth
quasi-compact $R$-model $X$ of $G$ such that the image of the
natural map $X(R)\rightarrow G(K)$ contains $<\!\!x\!\!>$, by
\cite[1.1.8+3.5.2]{neron}. By the universal property of the
N\'eron $lft$-model, there exists a unique $R$-morphism
$h:X\rightarrow \mathcal{G}^{lft}$ extending the isomorphism
$X_K\cong G$. Since $X$ is quasi-compact, the morphism $h_s$
factors through a finite union of connected components of
$\mathcal{G}^{lft}_s$. On the other hand, the image of
$h_s(k):X_s(k)\rightarrow \mathcal{G}_s^{lft}(k)$ contains the
subgroup of $\mathcal{G}^{lft}_s(k)$ generated by $\psi_s$, so the
connected component of $\mathcal{G}^{lft}_s$ containing $\psi_s$
is torsion in $\pi_0(\mathcal{G}^{lft}_s)$, and we have $\psi_s\in
\mathcal{G}_s(k)$ by Proposition \ref{prop-lftner}.
\end{proof}
In particular, if $K$ is complete, then $G^b(K')$ is the maximal
bounded subgroup of $G(K')$, for every finite unramified extension
$K'$ of $K$.
%

\begin{remark}
We do not know if the condition that $R$ is excellent is necessary
in Propositions \ref{prop-equiv} and \ref{prop-maxbound}. If
Proposition \ref{prop-unram} holds without the condition that $G$
admits a N\'eron $lft$-model, i.e., if N\'eron models always
commute with unramified base change, then the excellence condition
can be omitted in Propositions \ref{prop-equiv} and
\ref{prop-maxbound}, since it follows from \cite[10.2.2]{neron}
that $G$ admits a N\'eron $lft$-model iff $G\times_K \widehat{K}$
admits one.
\end{remark}


\section{Edixhoven's filtration and Chai's base change
conductor}\label{sec-edixchai} In this section, we assume that $K$
is strictly Henselian. In \cite{edix}, Edixhoven constructed a
filtration on the special fiber of the N\'eron model of an abelian
$K$-variety, which measures the behaviour of the N\'eron model
under finite tame extensions of $K$. This construction generalizes
without additional effort to the class of smooth and commutative
algebraic $K$-groups $G$ such that $G\times_K K'$ admits a N\'eron
model for all finite tame extensions $K'$ of $K$, and in
particular to the class of semi-abelian $K$-varieties. The
construction is explained in Section \ref{subsec-edix}, rephrased
in the language of Greenberg schemes. Most of the results in
Section \ref{subsec-edix} were stated (for abelian varieties) in
\cite{edix}, but many proofs were omitted. Since these results are
vital for the applications in this article, and since some of them
don't seem trivial to us (in particular Theorem \ref{theo-diff}),
we found it worthwhile to supply detailed proofs here if they were
not given in \cite{edix}. If $G$ is a tamely ramified semi-abelian
$K$-variety, then we relate Edixhoven's filtration to Chai's base
change conductor \cite{chai} in Sections \ref{subsec-chai} and
\ref{subsec-compar}.
\subsection{Edixhoven's filtration}\label{subsec-edix}
Let $K'$ be a tame finite extension of $K$ of degree $d$, and
denote by $R'$ the normalization of $R$ in $K'$. For each integer
$n> 0$, we put $R'_n=R'/(\mathfrak{M}')^{n}$, where
$\mathfrak{M}'$ is the maximal ideal of $R'$. We put $R'_0=\{0\}$,
the trivial $R'$-module. If $1\leq n\leq d$, then $R'_n$ carries a
natural $k$-algebra structure $k\cong R/\mathfrak{M}\rightarrow
R'/(\mathfrak{M}')^n$.

 Let $G$ be a smooth
commutative algebraic $K$-group such that $G'=G\times_K K'$ admits
a N\'eron model $\mathcal{G}'$.
 We denote by $X$ the Weil restriction $\prod_{R'/R}\mathcal{G}'$
of $\mathcal{G}'$ to $R$. By \cite[7.6.4]{neron} it is
representable by a group $R$-scheme, since $R'$ is finite and flat
over $R$ and $\mathcal{G}'$ is quasi-projective over $R'$
\cite[6.4.1]{neron}. The $R$-scheme $X$ is separated, smooth, and
of finite type over $R$ \cite[7.6.5]{neron}.

The extension $K'/K$ is Galois. We denote its Galois group by
$\mu$, and we let $\mu$ act on $K'$ from the left. The action of
$\zeta\in \mu$ on $\mathfrak{M}'/(\mathfrak{M}')^2$ is
multiplication by $\zeta'$, for some $\zeta'\in k$, and the map
$\zeta\mapsto \zeta'$ is an isomorphism between $\mu$ and the
group $\mu_d(k)$ of $d$-th roots of unity in $k$. By the universal
property of the N\'eron model, the right $\mu$-action on $G'$
extends uniquely to a $\mu$-action on the scheme $\mathcal{G}'$
such that the structural morphism $\mathcal{G}'\rightarrow \Spec
R'$ is $\mu$-equivariant. As in \cite[2.4]{edix}, this action
induces a right $\mu$-action on $X$.


The fixed point functor $(\cdot)^{\mu}$ from the category of
schemes with right $\mu$-action to the category of schemes is
right adjoint to the functor endowing a scheme with the trivial
$\mu$-action. We refer to \cite[3.1]{edix} for its basic
properties.

We have a tautological morphism of $K$-groups
$$G\rightarrow X_K\cong \prod_{K'/K}G'$$
By Galois descent, it factors through an isomorphism
$$G\rightarrow (X_K)^{\mu}\cong (X^{\mu})_K$$
\begin{prop}\label{prop-fix}
The $R$-scheme $X^{\mu}$ is a N\'eron model for $G$.
\end{prop}
\begin{proof}
 By \cite[3.4]{edix} the group
$R$-scheme $X^{\mu}$ is smooth. Let $H$ be any smooth group
$R$-scheme of finite type, and let
$$f:H_K\rightarrow (X^{\mu})_K$$ be a morphism of group $K$-schemes. Then
$f$ corresponds to a $\mu$-equivariant $K'$-morphism $f:H\times_R
K'\rightarrow G'$ where $\mu$ acts on $H\times_R K'$ via the
Galois action on $K'$. By the universal property of the N\'eron
model, $f$ extends to a $\mu$-equivariant morphism of group
$R'$-schemes $f':H\times_R R'\rightarrow \mathcal{G}'$, which
yields a morphism $g:H\rightarrow X$. If we let $\mu$ act
trivially on $H$, then $g$ is $\mu$-equivariant, so that $g$
factors through a morphism of group $R$-schemes $g:H\rightarrow
X^{\mu}$ extending $f$.
\end{proof}

\begin{definition}
For any $R'$-scheme $Y$ and any $i\in \{1,\ldots,d\}$ we put
$$Green_i(Y)=
\prod_{R'_{i}/k}(Y\times_{R'}R'_{i})$$ For $j\geq i$ in
$\{1,\ldots,d\}$ the truncation morphism $R'_j\rightarrow R'_i$
induces a morphism of $k$-schemes
$$\theta^j_i:Green_j(Y)\rightarrow
Green_i(Y)$$
\end{definition}
Since $\Spec R'_i\rightarrow \Spec k$ is universally bijective,
the proof of \cite[7.6.4]{neron} shows that $Green_i(Y)$ is indeed
representable by a $k$-scheme. The functor $Green_i(\cdot)$ is
compatible with open immersions. Note that the $k$-scheme
$Green_i(\mathcal{G}')$ inherits a group
 structure from $\mathcal{G}'$, as well as a right $\mu$-action.
 The truncation morphisms $\theta^j_i$ are $\mu$-equivariant
 morphisms of group $k$-schemes.

\begin{remark}
If $k$ is perfect or $R$ has equal characteristic, then by
\cite[4.1]{NiSe-weilres}, $Green_i(\mathcal{G}')$ is canonically
isomorphic to the Greenberg scheme $Gr^{R'}_{i-1}(\mathcal{G}')$
of $\mathcal{G}'$ (mind the shift of index in our notation $R'_i$
w.r.t. \cite{NiSe-weilres}).
\end{remark}

\begin{definition}
For $i\in \{1,\ldots,d\}$ we define $F^iX_s$ as the kernel of
$$\theta^d_i:X_s\cong Green_d(\mathcal{G}')\rightarrow Green_i(\mathcal{G}')$$
 This defines a decreasing filtration
$$X_s=:F^0X_s \supset F^1X_s\supset \ldots \supset F^dX_s=0$$
by subgroup $k$-schemes that are stable under the $\mu$-action.

For $i\in \{0,\ldots,d-1\}$, we denote by $Gr^iX_s$ the group
$k$-scheme $F^iX_s/F^{i+1}X_s$. It inherits a $\mu$-action from
$F^iX_s$.
\end{definition}

\begin{prop}\label{prop-triv}
Denote by $m$ the dimension of $G$. There exists a Zariski cover
$\mathcal{U}$ of $\mathcal{G}'$ such that for each member $U$ of
$\mathcal{U}$ and each pair of integers $j\geq i$ in
$\{1,\ldots,d\}$, the truncation morphism
$$\theta^{j}_i:Green_{j}(U)\rightarrow
Green_i(U)$$ is a trivial fibration whose fiber is isomorphic (as
a $k$-scheme) to $\A_k^{(j-i)m}$.

For each $i\in \{1,\ldots,d-1\}$, the kernel of $\theta^{i+1}_i$
is canonically isomorphic to $Gr^iX_s$. Moreover, the group
$k$-scheme
 $Gr^0X_s$ is canonically isomorphic to $\mathcal{G}'_s$.
\end{prop}
\begin{proof}
The fact that $\theta^j_i$ is a trivial fibration with fiber
$\A_k^{(j-i)m}$ above the members of a Zariski cover of
$\mathcal{G}'$ can be proven exactly as in the case where $k$ is
perfect \cite[3.4.2]{sebag1}. For the second part of the
statement, consider the short exact sequence of group $k$-schemes
$$\begin{CD}0@>>>
F^iX_s@>>> X_s=Green_d(\mathcal{G}')@>\theta^d_i>>
Green_i(\mathcal{G}')@>>> 0\end{CD}$$ Dividing by $F^{i+1}X_s$ we
obtain an exact sequence
$$\begin{CD}0@>>>
Gr^iX_s@>>> X_s/F^{i+1}X_s @>\overline{\theta}^d_i>>
Green_i(\mathcal{G}')@>>> 0\end{CD}$$ However, by the first exact
sequence (with $i$ replaced by $i+1$) we see that there exists a
canonical isomorphism $X_s/F^{i+1}X_s\cong
Green_{i+1}(\mathcal{G}')$ identifying $\overline{\theta}^d_i$
with $\theta^{i+1}_i$. Considering the short exact sequence
$$\begin{CD}0@>>> F^1X_s@>>> X_s@>\theta^d_1>> Green_1(\mathcal{G}')\cong
\mathcal{G}'_s@>>> 0\end{CD}$$ we see that $Gr^0X_s\cong
\mathcal{G}'_s$.
\end{proof}

\begin{prop}\label{prop-quot}
For each $0<i<d$, there exists a canonical $\mu$-equivariant
isomorphism of group $k$-schemes
$$Gr^iX_s\cong Lie(\mathcal{G}'_s)\otimes_k
(\mathfrak{M}'/(\mathfrak{M}')^2)^{\otimes i}$$ where we view the
right hand side as a vector group $k$-scheme, and where the right
action of $\mu$ on $(\mathfrak{M}'/(\mathfrak{M}')^2)^{\otimes i}$
is the inverse of the left Galois action. In particular, $F^iX_s$
is unipotent, smooth and connected for $0<i< d$.
\end{prop}
\begin{proof}
The proof of \cite[\S 5.1]{edix} carries over without changes. If
$k$ is perfect or $R$ has equal characteristic, see also \cite[\S
2]{greenbergII}.
\end{proof}

\begin{definition}
We define a decreasing filtration on $\mathcal{G}_s$ by subgroup
$k$-schemes
$$\mathcal{G}_s=F^{0}_{d}\mathcal{G}_s\supset F^1_d\mathcal{G}_s\supset \ldots \supset F^d_d\mathcal{G}_s=0$$
 by putting
$F^i_d\mathcal{G}_s=(F^iX_s)^{\mu}$. For each $i\in
\{0,\ldots,d-1\}$, we put
$$Gr^i_d\mathcal{G}_s=F^i_d\mathcal{G}_s/F^{i+1}_d\mathcal{G}_s$$
\end{definition}
Note that we indeed have
$\mathcal{G}_s=(F^0X_s)^{\mu}=(X_s)^{\mu}$ by Proposition
\ref{prop-fix}.
 It follows immediately from the definition that
$F^i_d\mathcal{G}_s$ is the kernel of the truncation morphism
$$\mathcal{G}_s\cong (X_s)^{\mu}\rightarrow
Green_i(\mathcal{G}')$$ for $i=1,\ldots,d$. Observe that, since
the extension $K'$ of $K$ of degree $d\in \N'$ is uniquely
determined by $d$ up to $K$-automorphism, the filtration
$F^\bullet_d\mathcal{G}_s$ only depends on $d$ and not on $K'$.

\begin{lemma}\label{lemm-exact}
For each $i\in \{0,\ldots,d-1\}$ there is a canonical isomorphism
$$Gr^i_d\mathcal{G}_s\cong (Gr^i X_s)^{\mu}$$
\end{lemma}
\begin{proof}
It suffices to show that the exact sequence
$$0\rightarrow F^{i+1}X_s\rightarrow F^iX_s\rightarrow
Gr^iX_s\rightarrow 0$$ remains exact after applying the fixed
point functor $(\cdot)^{\mu}$, for each $i\in \{0,\ldots,d-1\}$.

Left exactness is clear. It remains to show that
$(F^iX_s)^{\mu}\rightarrow (Gr^iX_s)^{\mu}$ is a surjection of
$fpqc$ sheaves. For any commutative group $k$-scheme $H$, we
denote by $d_H:H\rightarrow H$ the multiplication by $d$. By
Proposition \ref{prop-quot}, $F^{i+1}X_s$ is unipotent, so since
$d$ is invertible in $k$, multiplication by $d$ is an automorphism
on $F^{i+1}X_s$. We denote its inverse by $(d^{-1})_{F^{i+1}X_s}$.

For any commutative group $k$-scheme $Z$ endowed with a right
$\mu$-action, consider the morphism $N_Z:Z\rightarrow Z^{\mu}$
defined by
$$N_Z(S):Z(S)\rightarrow Z^{\mu}(S):s\mapsto \sum_{\zeta\in \mu}s*\zeta $$ for all $k$-schemes $S$.
If we denote by $\iota_Z$ the tautological closed immersion
$Z^{\mu}\rightarrow Z$, then  $N_Z\circ \iota_Z=d_{Z^{\mu}}$.

Let $S$ be a $k$-scheme, and $c$ a section in
$(Gr^iX_s)^{\mu}(S)$.
 Choose an $fpqc$ covering $S'\rightarrow S$ such that $c$ lifts to
 an element $b$ of $(F^iX_s)(S')$. Since $c$ is invariant
 under the $\mu$-action, the element $$b'=N_{F^iX_s}(b)-d_{F^i X_s}(b)$$ maps to zero in
 $(Gr^iX_s)(S')$, so it belongs to $(F^{i+1}X_s)(S')$. Now
 $$b''=b+(d^{-1})_{F^{i+1}X_s}(b')$$ maps to $c$. This element belongs to
 $(F^iX_s)^{\mu}(S')$,
 since for any $\zeta\in \mu$, $b''* \zeta-b''$ is an element
 of $(F^{i+1}X_s)(S')$ that is mapped to zero by $d_{F^{i+1}X_s}$, so that $b''* \zeta=b''$.
\end{proof}
\begin{cor}\label{cor-fixed}
 The group $k$-scheme $Gr^0_d \mathcal{G}_s$ is canonically isomorphic to
 $(\mathcal{G}'_s)^{\mu}$. It is a smooth algebraic $k$-group, and
 there is a canonical isomorphism of $k$-vector spaces
 $$Lie((\mathcal{G}'_s)^{\mu})\cong (Lie(\mathcal{G}'_s))^{\mu}$$
 For each
$0<i<d$, there exists a canonical isomorphism
\begin{equation}\label{eq-twist}
Gr^i_d\mathcal{G}_s\cong Lie(\mathcal{G}'_s)[i]\otimes_k
(\mathfrak{M}'/(\mathfrak{M}')^2)^{\otimes i}\end{equation} where
$Lie(\mathcal{G}'_s)[i]$ is the subspace of $Lie(\mathcal{G}'_s)$
where $\mu\cong \mu_d(k)$ acts as $(v,\zeta)\mapsto \zeta^i\cdot
v$ for $(v,\zeta)\in  Lie(\mathcal{G}'_s)\times \mu_d(k)$, and
where we view the right hand side of (\ref{eq-twist}) as a vector
group $k$-scheme. In particular, $F^i_d\mathcal{G}_s$ is
unipotent, smooth and connected for $0<i\leq d$.
\end{cor}
\begin{proof}
Smoothness of $(\mathcal{G}'_s)^{\mu}$ and the isomorphism
 $$Lie((\mathcal{G}'_s)^{\mu})\cong (Lie(\mathcal{G}'_s))^{\mu}$$
 follow from \cite[3.2+4]{edix}. The remaining statements follow immediately from
 Proposition \ref{prop-quot} and Lemma \ref{lemm-exact}.
\end{proof}

\begin{definition}
We say that an integer $j$ in $\{0,\ldots,d-1\}$ is a $K'$-jump of
$G$ if $Gr^j_d\mathcal{G}_s\neq 0$. The dimension of
$Gr^j_d\mathcal{G}_s$ is called the multiplicity of the jump $j$.
\end{definition}

By Corollary \ref{cor-fixed}, the $K'$-jumps of $G$ and their
multiplicites can be computed from the $\mu$-action on
$Lie(\mathcal{G}'_s)$.


The following theorem will play a crucial role in the remainder of
this article. If $G$ is an abelian variety, the result was stated
in \cite[5.4.6]{edix} without proof.
\begin{theorem}\label{theo-diff}
Let $K'$ be a finite tame extension of $K$, and denote by $R'$ the
normalization of $R$ in $K'$. Let $G$ be a smooth commutative
algebraic $K$-group, and assume that $G'=G\times_K K'$ admits a
N\'eron model $\mathcal{G}'$.  Denote by $\mathcal{G}$ the N\'eron
model of $G$, and by $\mathcal{K}$ the kernel of the natural
morphism $h:\mathcal{G}\times_R R'\rightarrow \mathcal{G}'$. Let
$j_1(G,K'),\ldots,j_u(G,K')$ be the $K'$-jumps associated to $G$,
with respective multiplicities $m_1(G,K'),\ldots,m_u(G,K')$.

Then the pull-back of the fundamental exact sequence of
$\mathcal{O}_{\mathcal{G}\times_R R'}$-modules
$$ h^*\Omega^1_{\mathcal{G}'/R'}\rightarrow
\Omega^1_{\mathcal{G}\times_R R'/R'}\rightarrow
\Omega^1_{\mathcal{G}\times_R R'/\mathcal{G}'}\rightarrow 0$$
w.r.t. the unit section $e_{\mathcal{G}\times_R R'}$ yields a
short exact sequence of $R'$-modules
\begin{equation}\label{seq}0\rightarrow
\omega_{\mathcal{G}'/R'}\rightarrow \omega_{\mathcal{G}\times_R
R'/R'}\rightarrow \omega_{\mathcal{K}/R'}\rightarrow
0\end{equation} and $\omega_{\mathcal{K}/R'}$ is isomorphic to
$$\left(\bigoplus_{m_1(G,K')}R'_{j_1(G,K')}\right)\oplus \cdots
\oplus \left(\bigoplus_{m_u(G,K')}R'_{j_u(G,K')}\right)$$
Moreover, $Lie(h)$ injects $Lie(\mathcal{G}\times_R R')$ into
$Lie(\mathcal{G}')$ and there exists an isomorphism of
$R'$-modules
\begin{equation}\label{eq-iso}\omega_{\mathcal{K}/R'}\cong
Lie(\mathcal{G}')/Lie(\mathcal{G}\times_R R')\end{equation}
\end{theorem}
\begin{proof}
We refer to \cite[\S\,1]{liu-lorenzini-raynaud} for some basic
results on Lie algebras of group schemes. To show that the
sequence (\ref{seq}) is exact, it suffices to show that
$$\omega_{\mathcal{G}'/R'}\rightarrow \omega_{\mathcal{G}\times_R
R'/R'}$$ is injective. This follows immediately from the fact that
$\omega_{\mathcal{G}'/R'}$ is free and $h_{K'}$ an isomorphism.
Dualizing (\ref{seq}) we find an exact sequence
$$\begin{CD}0@>>> Lie(\mathcal{G}\times_R R')@>Lie(h)>>
Lie(\mathcal{G}')@>>> Ext^1_{R'}(\omega_{\mathcal{K}/R'},R')@>>>
0\end{CD}$$ Since $\omega_{\mathcal{K}/R'}$ is torsion, it is
easily seen that the $R'$-module
$Ext^1_{R'}(\omega_{\mathcal{K}/R'},R')$ is isomorphic to
$\omega_{\mathcal{K}/R'}$.

We put $d=[K':K]$. We denote by $\mathfrak{M}'$ the maximal ideal
of $R'$, and by $\mu$ the Galois group of the extension $K'/K$. We
let $\mu$ act on $K'$ from the left.

We put $X=\prod_{R'/R}\mathcal{G}'$ as before. Consider the
commutative diagram
$$\begin{CD}
\mathcal{G}@>\alpha >> X^{\mu}
\\ @V\beta VV @VV\gamma V
\\ \prod_{R'/R}(\mathcal{G}\times_R R')@>>\delta=\prod_{R'/R}h> X
\end{CD}$$
where $\alpha$ is the isomorphism from Proposition \ref{prop-fix},
and
$\beta$ and $\gamma$ are the tautological morphisms. 
 It is easily seen that $Lie(\cdot)$ commutes with Weil
restriction, so we have canonical isomorphisms of $R$-modules
\begin{eqnarray*}
 Lie(\prod_{R'/R}(\mathcal{G}\times_R R'))&\cong &
Lie(\mathcal{G}\times_R R')
\\ Lie(X)&\cong& Lie(\mathcal{G}')
\end{eqnarray*}
that identify $\delta$ with $Lie(h)$.

Arguing as in the proof of \cite[3.2]{edix} we see that
$Lie(\gamma)$ is an isomorphism onto $Lie(\mathcal{G}')^{\mu}$, so
that we can use $Lie(\gamma)\circ Lie(\alpha)$ to identify
$Lie(\mathcal{G})$ with $Lie(\mathcal{G}')^{\mu}$. Also, taking
formal parameters, it is easily seen that the morphism of
$R'$-modules
$$Lie(\mathcal{G}')^{\mu}\otimes_R R'\rightarrow Lie(\mathcal{G}\times_R
R')$$ obtained from $Lie(\beta)$ by extension of scalars is an
isomorphism. Modulo this isomorphism, the morphism of $R'$-modules
$Lie(h):Lie(\mathcal{G}\times_R R')\rightarrow Lie(\mathcal{G}')$
corresponds to the inclusion
$$Lie(\mathcal{G}')^{\mu}\otimes_R R'\rightarrow
Lie(\mathcal{G}')$$

Let $\{\overline{e}_1,\ldots,\overline{e}_m\}$ be a $k$-basis of
$Lie(\mathcal{G}'_s)$ such that $\mu\cong \mu_d(k)$ acts on
$\overline{e}_i$ by $$\overline{e}_i*\zeta=\zeta^{a_i}\cdot
\overline{e}_i$$ for each $\zeta\in \mu_d(k)$ and each
$i\in\{1,\ldots,m\}$, with $a_i\in \{0,\ldots,d-1\}$. Then by
Corollary \ref{cor-fixed}, the occurring $a_i$ (with
multiplicities) are the $K'$-jumps of $G$ (with multiplicities),
so we have to show that there exists an isomorphism
$$Lie(\mathcal{G}')/(Lie(\mathcal{G}')^{\mu}\otimes_R R')\cong \oplus_{i=1}^{m}R'_{a_i}$$

We have a natural isomorphism of $k$-vector spaces
$$Lie(X_s)\cong Lie(\mathcal{G}'\times_{R'}R'_d)$$
Consider the projection
$$Lie(\theta^d_0):Lie(X_s)\rightarrow Lie(\mathcal{G}'_s)$$ and
 lift the basis
$\{\overline{e}_1,\ldots,\overline{e}_m\}$ of
$Lie(\mathcal{G}'_s)$ to a tuple $\{e_1,\ldots,e_m\}$ of elements
in $Lie(X_s)$ such that $\mu\cong \mu_d(k)$ acts on $e_i$ by
$$e_i*\zeta=\zeta^{a_i}\cdot e_i$$ for each $\zeta\in \mu_d(k)$
and each $i\in \{1,\ldots,m\}$. Choose a uniformizer $\pi'$ in
$R'$. Since $\{e_1,\ldots,e_m\}$ is a $R'_d$-basis for
$Lie(\mathcal{G}'\times_{R'}R'_d)$, we see that
$$\{v_{i,j}=e_i\cdot (\pi')^j\,|\,i=1,\ldots,m,\ j=0,\ldots,d-1\}$$ is a
$k$-basis of $Lie(\mathcal{G}'\times_{R'}R'_d)\cong Lie(X_s)$, and
$\zeta\in \mu_d(k)$ acts on $v_{i,j}$ by
$$v_{i,j}*\zeta=\zeta^{a_i-j}\cdot v_{i,j}$$
(see the construction of the group action in \cite[2.4]{edix} for
the origin of the sign in the exponent $-j$).

If $V_1$ and $V_2$ are (right) $R[\mu]$-modules of finite type,
free over $R$, and $\varphi:V_1\otimes_R k\rightarrow V_2\otimes_R
k$ is an isomorphism of  $k[\mu]$-modules, then $\varphi$ lifts to
an isomorphism of $R[\mu]$-modules $\varphi':V_1\rightarrow V_2$.
To see this, take any morphism of $R$-modules $\phi:V_1\rightarrow
V_2$ lifting $\varphi$. The morphism
$$\varphi':V_1\rightarrow V_2:v\mapsto \frac{1}{d}\sum_{\zeta\in
\mu}\phi(v*\zeta)*\zeta^{-1}$$ is a morphism of $R[\mu]$-modules.
It is an isomorphism since its reduction modulo $\pi$ is the
isomorphism $\varphi$.

 This implies that we can lift the $k$-basis $\{v_{i,j}\}$ of $Lie(X_s)$ to a $R$-basis $\{w_{i,j}\}$
of $Lie(\mathcal{G}')$ such that the action of $\mu\cong
\mu_d(k)\cong \mu_d(R)$ on $w_{i,j}$ is given by
$$ w_{i,j}*\xi=\xi^{a_i-j}\cdot w_{i,j}$$ for each $\xi\in \mu_d(R)$ and all $i,j$. We see
that the $R$-module $$Lie(\mathcal{G})=Lie(\mathcal{G'})^{\mu}$$
is generated by the elements $w_{i,a_i}$ with $i=1,\ldots,m$. We
observe, moreover, that $$\pi \cdot Lie(\mathcal{G}')\subset
Lie(\mathcal{G}')^{\mu}\otimes_R R'$$ because
$Lie(\mathcal{G}')^{\mu}\otimes_R R'$ is a sub-$R'$-module of
$Lie(\mathcal{G}')$ of the same rank, and, for any $v$ in
$Lie(\mathcal{G}')$, we have $\pi\cdot v\in
Lie(\mathcal{G}')^{\mu}$ iff $v\in Lie(\mathcal{G}')^{\mu}$.

 Hence, we can conclude that
\begin{eqnarray*}
Lie(\mathcal{G}')/(Lie(\mathcal{G}')^{\mu}\otimes_R R')&\cong&
Lie(\mathcal{G}'\times_{R'} R'_d)/(\sum_{i=1}^{m} R'_d\cdot
v_{i,a_i})
\\ &\cong& Lie(\mathcal{G}'\times_{R'} R'_d)/(\sum_{i=1}^{m} (\pi')^{a_i} R'_d\cdot  e_{i})
\\ &\cong& \oplus_{i=1}^{m}R'_{a_i}
\end{eqnarray*}
\end{proof}

\begin{lemma}\label{lemm-indep}
Let $G$ be a smooth commutative algebraic $K$-group, and assume
that $G\times_K K'$ admits a N\'eron model for every finite tame
extension $K'$ of $K$. For all $d,n\in \N'$ and each $i\in
\{0,\ldots,d\}$ we have
$F^i_d\mathcal{G}_s=F^{in}_{dn}\mathcal{G}_s$.
\end{lemma}
\begin{proof}
For $i=0$ the statement is obvious, so assume $i>0$. Let
$K'\subset K''$ be extensions of $K$ of degree $d$, resp. $dn$,
with rings of integers $R'$, resp., $R''$. Denote  by
$\mathcal{G}'$ and $\mathcal{G}''$ the N\'eron models of
$G\times_K K'$, resp. $G\times_K K''$. It suffices to show that
the kernels of the truncation morphisms
\begin{eqnarray*}
\phi_1&:&\mathcal{G}_s\rightarrow Green_{i}(\mathcal{G}')
\\ \phi_2&:&\mathcal{G}_s\rightarrow
Green_{in}(\mathcal{G}'')
\end{eqnarray*}
coincide. We have a natural morphism
$\mathcal{G}'\times_{R'}R''\rightarrow \mathcal{G}''$ inducing a
morphism $\mathcal{G}'\rightarrow \prod_{R''/R'}\mathcal{G}''$,
which is a closed immersion by Proposition \ref{prop-fix}.
Moreover, $Green_{in}(\mathcal{G}'')$ is canonically isomorphic to
$Green_{i}(\prod_{R''/R'}\mathcal{G}'')$, and since Weil
restriction respects closed immersions, we obtain a closed
immersion of group $k$-schemes
$$\phi_{1,2}:Green_{i}(\mathcal{G}')\rightarrow Green_{in}(\mathcal{G}'')$$
such that $\phi_2=\phi_{1,2}\circ\phi_1$.
%
\end{proof}

\begin{definition}
Let $G$ be a smooth commutative algebraic $K$-group, and assume
that $G\times_K K'$ admits a N\'eron model for every finite tame
extension of $K$. Denote by $q$ the characteristic of $k$. For
each element $\alpha=a/b$ of $\Z_{(q)}\cap[0,1[$, with $a\in \N$
and $b\in \N'$, we put
$\widetilde{F}^\alpha\mathcal{G}_s=F^{a}_b\mathcal{G}_s$. By Lemma
\ref{lemm-indep}, this definition does not depend on the choice of
$a$ and $b$. Then $\widetilde{F}^{\bullet}\mathcal{G}_s$ is a
decreasing filtration on $\mathcal{G}_s$ by subgroup $k$-schemes.

Let $\rho$ be an element of $\R\cap[0,1[$. We put
$\widetilde{F}^{>\rho}\mathcal{G}_s=\widetilde{F}^{\beta}\mathcal{G}_s$,
where $\beta$ is any value in $\Z_{(q)}\cap \,]\rho,1[$ such that
$\widetilde{F}^{\beta'}\mathcal{G}_s=\widetilde{F}^{\beta}\mathcal{G}_s$
for all $\beta'$ in $\Z_{(q)}\cap \,]\rho,\beta]$. If $\rho\neq 0$
we put
$\widetilde{F}^{<\rho}\mathcal{G}_s=\widetilde{F}^{\gamma}\mathcal{G}_s$
where $\gamma$ is any value in $\Z_{(q)}\cap [0,\rho[$ such that
$\widetilde{F}^{\gamma'}\mathcal{G}_s=\widetilde{F}^{\gamma}\mathcal{G}_s$
for all $\gamma'$ in $\Z_{(q)}\cap [\gamma,\rho[$. We put
$\widetilde{F}^{<0}\mathcal{G}_s=\mathcal{G}_s$.

 We define
$$\widetilde{Gr}^{\rho}\mathcal{G}_s=\widetilde{F}^{<\rho}\mathcal{G}_s/\widetilde{F}^{>\rho}\mathcal{G}_s$$
We say that $j\in \R\cap [0,1[$ is a jump of $G$ if
$\widetilde{Gr}^{j}\mathcal{G}_s\neq 0$. The multiplicity of $j$
is the dimension of $\widetilde{Gr}^{j}\mathcal{G}_s$. We say that
the multiplicity of $j$ as a jump of $G$ is zero if $j$ is not a
jump of $G$.
\end{definition}
It follows immediately from the definition that the sum of the
multiplicities of the jumps of $G$ equals the dimension of $G$.
%
As noted by Edixhoven in \cite[5.4.5]{edix}, it is not clear if
the jumps of $G$ are rational numbers, but one can be more precise
if $G$ is a tamely ramified semi-abelian variety over $K$.

\begin{prop}\label{prop-tamejump}
Denote by $q$ the characteristic of $k$. Assume that $G$ is a
semi-abelian $K$-variety, and that $G$ acquires semi-abelian
reduction over a tame finite extension $K'$ of $K$ of degree $d$.
For any $\alpha\in \Z_{(q)}\cap [0,1[$ we have
$$\widetilde{F}^{<\alpha}\mathcal{G}_s=\widetilde{F}^{\alpha}\mathcal{G}_s=F^{\lceil \alpha\cdot d\rceil}_d\mathcal{G}_s$$

 If we denote by
$j_1(G,K'),\ldots,j_u(G,K')$ the $K'$-jumps of $G$, with
multiplicities $m_1(G,K'),\ldots,m_u(G,K')$, then the jumps of $G$
are given by $$j_1(G,K')/d,\ldots,j_u(G,K')/d$$ with the same
multiplicities. In particular, the jumps of $G$ belong to
$\Z[1/d]\cap [0,1[$.

Moreover, if $L/K$ is any finite tame extension of $K$, of degree
$e$, then the set of $L$-jumps of $G$ is
$$\{\lfloor j_1(G,K')\cdot e/d \rfloor ,\ldots,\lfloor j_u(G,K')\cdot e/d\rfloor  \}$$
 The multiplicity of $j\in \{0,\ldots,e-1\}$
as a $L$-jump of $G$ equals
$$\sum_{i\in T_j}m_i(G,K')$$
with $T_j=\{i\in\{1,\ldots,u\}\,|\,j=\lfloor j_i(G,K')\cdot e/d
\rfloor\}$.
\end{prop}
\begin{proof}
Write $\alpha$ as $i/(dn)$ with $d,n\in \N'$. We'll show that
\begin{equation}\label{eq-tamej}\widetilde{F}^{\alpha}\mathcal{G}_s:=F^{i}_{dn}\mathcal{G}_s=F^{\lceil
i/n\rceil}_{d}\mathcal{G}_s\end{equation}
The remainder of the statement follows
 from (\ref{eq-tamej}) by some elementary combinatorics.

So let us prove (\ref{eq-tamej}). Let $K'\subset K''$ be an
extension of degree $n$. Denote by $R'$ and $R''$ the ring of
integers of $K'$, resp. $K''$, and by $\mathcal{G}'$ and
$\mathcal{G}''$ the N\'eron models of $G\times_K K'$, resp.
$G\times_K K''$. Since $G$ has semi-abelian reduction over $K'$,
the natural morphism $\mathcal{G}'\times_{R'}R''\rightarrow
\mathcal{G}''$ is an open immersion, and $Lie(\mathcal{G}'_s)$ and
$Lie(\mathcal{G}''_s)$ are isomorphic as $k$-vector spaces with
$\mu''=G(K''/K)$-action. In particular, $G(K''/K')$ acts trivially
on $Lie(\mathcal{G}''_s)$. By Corollary \ref{cor-fixed}, this
means that $Gr^j_{nd}\mathcal{G}_s=0$ if $j$ is not divisible by
$n$. It follows that
$$F^{j}_{dn}\mathcal{G}_s=F^{\lceil j/n\rceil\cdot
n}_{dn}\mathcal{G}_s=F^{\lceil j/n\rceil}_{d}\mathcal{G}_s$$ for
$j\in \{0,\ldots,dn\}$.
\end{proof}
\begin{remark}
If $G$ is the Jacobian of a smooth and proper $K$-curve with a
$K$-rational point, then the jumps of $G$ are rational, without
any tameness condition. In fact, much more can be said: see
\cite[8.4]{halle-neron}.
\end{remark}

\subsection{Chai's base change conductor}\label{subsec-chai}
Let $G$ be a smooth commutative algebraic $K$-group that admits a
N\'eron model $\mathcal{G}$. Let $K'$ be a finite separable
extension of $K$ of degree $d$, and denote by $R'$ the
normalization of $R$ in $K'$. Assume that $G'=G\times_K K'$ admits
a N\'eron model $\mathcal{G}'$. By the universal property of the
N\'eron model, there exists a unique morphism of group
$R'$-schemes
$$h:\mathcal{G}\times_R R'\rightarrow \mathcal{G}'$$
that extends the canonical isomorphism between the generic fibers.
Since $h_{K'}$ is an isomorphism, the map $Lie(h)$ injects
$Lie(\mathcal{G}\times_R R')$ into $Lie(\mathcal{G}')$, and the
quotient $R'$-module
$$Lie(\mathcal{G}')/(Lie(\mathcal{G}\times_R R'))$$
is a torsion module of finite type over $R'$. If we denote by
$\mathfrak{M}'$ the maximal ideal of $R'$ and if we put
$R'_i=R'/(\mathfrak{M}')^i$ for each $i\in \Z_{>0}$, we get a
decomposition
$$Lie(\mathcal{G}')/(Lie(\mathcal{G}\times_R R'))\cong \bigoplus_{i=1}^v R'_{ c_i(G,K')\cdot d}$$
with $0<c_1(G,K')\leq \ldots\leq c_v(G,K')$ in $\Z[1/d]$.

\begin{definition}
We call the tuple of rational numbers
$(c_1(G,K'),\ldots,c_v(G,K'))$ the tuple of $K'$-elementary
divisors associated to $G$, and we call
$$c(G,K')=c_1(G,K')+\cdots +c_v(G,K')=\frac{1}{d}\cdot \mathrm{length}_{R'}\left(Lie(\mathcal{G}')/(Lie(\mathcal{G}\times_R R'))\right)$$ the $K'$-base change conductor
associated to $G$.
\end{definition}

As a special case, we recall the following definition from
\cite[2.4]{chai}.

\begin{definition}[Chai]
Let $A$ be a semi-abelian variety, and let $K'$ be a finite
separable extension of $K$ such that $A\times_K K'$ has
semi-abelian reduction. The values $c_i(A,K')$ and $c(A,K')$ only
depend on $A$, and not on $K'$. We call them the elementary
divisors, resp. the base change conductor of $A$.
\end{definition}
Note that our definition differs slightly from the one in
\cite[2.4]{chai}. Chai extends the tuple of elementary divisors by
adding zeroes to the left until the length of the tuple equals the
dimension of $A$. Our definition is more convenient for the
purpose of this paper.

\begin{prop}\label{prop-c0}
For any semi-abelian variety $A$, we have $c(A)=0$ iff $A$ has
semi-abelian reduction.
\end{prop}
\begin{proof}
The ``if'' part is obvious, so let us prove the converse
implication. Take a finite separable extension $K'$ of $K$ such
that $A'=A\times_K K'$ has semi-abelian reduction. Denote by $R'$
the normalization of $R$ in $K'$ and by $\mathcal{A}$ and
$\mathcal{A}'$ the N\'eron models of $A$, resp. $A'$. Consider the
canonical morphism $h:\mathcal{A}\times_R R'\rightarrow
\mathcal{A}'$. Then $c(A)=0$ implies that $Lie(h)$ is an
isomorphism. Hence, $h$ is \'etale, and since $h_K$ is an
isomorphism, $h$ is an open immersion \cite[2.3.2']{neron}.
Therefore, $A$ has semi-abelian reduction.
\end{proof}
\subsection{A comparison result}\label{subsec-compar}
\begin{theorem}\label{theo-compar}
Let $K'$ be a finite tame extension of $K$, of degree $d$. Let $G$
be a smooth commutative algebraic $K$-group, and assume that $G$
and  $G'=G\times_K K'$ admit N\'eron models.  Let
$j_1(G,K')<\ldots <j_u(G,K')$ be the non-zero $K'$-jumps of $G$,
with respective multiplicities $m_1(G,K'),\ldots,m_u(G,K')$.

Then the tuple of  $K'$-elementary divisors $c_i(G,K')$ associated
to $G$ is
$$\left(\underbrace{j_1(G,K')/d,\ldots,j_1(G,K')/d}_{m_1(G,K')\times},
\ldots,
\underbrace{j_u(G,K')/d,\ldots,j_u(G,K')/d}_{m_u(G,K')\times}\right)$$
and we have
$$d\cdot c(G,K')=\sum_{i=1}^{u}(m_i(G,K')\cdot j_i(G,K'))$$
\end{theorem}
\begin{proof}
This follows immediately from Theorem \ref{theo-diff}.
\end{proof}
Combined with the fact that the sum of the multiplicities of the
$K'$-jumps of $G$ equals the dimension of $G$, Theorem
\ref{theo-compar} shows that the $K'$-jumps (with multiplicities)
and the $K'$-elementary divisors determine each other.
\begin{cor}\label{cor-compar}
Let $A$ be a tamely ramified semi-abelian $K$-variety, and let
$j_1(A)<\ldots<j_u(A)$ be the non-zero jumps associated to $A$,
with respective multiplicities $m_1(A),\ldots,m_u(A)$. Then the
tuple of elementary divisors of $A$ is given by
$$\left(\underbrace{j_1(A),\ldots,j_1(A)}_{m_1(A)\times},\ldots,\underbrace{j_i(A),\ldots,j_i(A)}_{m_i(A)\times},
\ldots, \underbrace{j_u(A),\ldots,j_u(A)}_{m_u(A)\times}\right)$$
 and we have
$$c(A)=\sum_{i=1}^{u}(m_i(A)\cdot j_i(A))$$
\end{cor}
\begin{proof}
Apply Theorem \ref{theo-compar} and Proposition
\ref{prop-tamejump}.
\end{proof}
\begin{cor}
If $A$ is a tamely ramified semi-abelian $K$-variety, then its
elementary divisors are contained in the interval $]0,1[$.
\end{cor}
\begin{cor}\label{cor-jumps0}
If $A$ is a tamely ramified semi-abelian $K$-variety, then $A$ has
semi-abelian reduction iff $0$ is the only jump of $A$.
\end{cor}
\begin{cor}
Let $A$ be a tamely ramified semi-abelian $K$-variety, and let
$K'$ be a finite tame extension of $K$ such that $A$ acquires
semi-abelian reduction over $K'$. Put $e=[K':K]$ and
$\mu=G(K'/K)$, and denote by $\mathcal{A}'$ the N\'eron model of
$A\times_K K'$. Then for any $\zeta\in \mu\cong \mu_e(k)$, the
determinant of the action of $\zeta$ on $Lie(\mathcal{A}'_s)$
equals $\zeta^{e\cdot c(A)}$.
\end{cor}
\begin{proof}
This follows from Corollary \ref{cor-fixed}, Proposition
\ref{prop-tamejump} and Corollary \ref{cor-compar}.
\end{proof}


\section{Jumps and monodromy eigenvalues}\label{sec-monodromy}
In this section, we assume that $K$ is strictly henselian and that
$k$ is algebraically closed. Recall that for any real number $x$
we denote by $[x]$ its decimal part $[x]=x-\lfloor x\rfloor\in
[0,1[$.

The following lemma describes the behaviour of the jumps of a
tamely ramified semi-abelian $K$-variety under tame base change.
For later use, we remark that the proof remains valid if $K$ is
strictly henselian but $k$ is not necessarily perfect.
\begin{lemma}\label{lemma-basech}
Let $A$ be a semi-abelian $K$-variety, and assume that $A$
acquires semi-abelian reduction on some finite tame extension $L$
of $K$, of degree $e$. Let $K'$ be a finite tame extension of $K$,
of degree $d$. Let $j_1(A),\ldots,j_u(A)$ be the jumps of $A$,
with respective multiplicities $m_1(A),\ldots,m_u(A)$. Then the
set of jumps of $A'=A\times_K K'$ is
$$J=\{[d\cdot j_1(A)],\ldots,[d\cdot j_u(A)]\}\subset \Z[1/e]\cap [0,1[$$
Moreover, for any $j\in J$, the multiplicity of $j$ as a jump of
$A'$ equals
$$\sum_{i\in S_j}m_i(A)$$
where $S_j$ is the subset of $\{1,\ldots,u\}$ consisting of
indices $i$ such that $[d\cdot j_i(A)]=j$.
\end{lemma}
\begin{proof}
We may assume that $K'$ is contained in $L$. Put $n=e/d$. Denote
by $\mathcal{B}$ the N\'eron model of $B=A\times_K L$, and put
$\mu=G(L/K)\cong \mu_{e}(k)$ and $\mu'=G(L/K')\cong \mu_{n}(k)$.
For $a=0,\ldots,e-1$ we denote by $V[a]$ the subspace of
$Lie(\mathcal{B}_s)$ where each $\zeta\in \mu$ acts by
$v*\zeta=\zeta^a\cdot v$. Likewise, for $b=0,\ldots,n-1$ we denote
by $V'[b]$ the subspace of $Lie(\mathcal{B}_s)$ where each $\xi\in
\mu'$ acts by $v*\xi=\xi^b\cdot v$. Then we obviously have
$$V'[b]=\bigoplus_{a\equiv b\,\mathrm{mod}\,n}V[a]$$
By Proposition \ref{prop-tamejump} and Corollary \ref{cor-fixed}
we know that $a/e$ is a jump of $A$ with multiplicity $m>0$ iff
$V[a]$ has dimension $m$, and likewise, $b/n$ is a jump of $A'$
with multiplicity $m'>0$ iff $V'[b]$ has dimension $m'$, so the
result follows.
\end{proof}

\begin{lemma}\label{lemma-subsemiab}
Let $F$ be a perfect field, and let $A$ be a semi-abelian
$F$-variety. If $B$ is a connected smooth subgroup $F$-scheme of
$A$, then $B$ is semi-abelian.
\end{lemma}
\begin{proof}
By \cite[2.3]{conrad-chevalley} there exists a unique connected
smooth linear subgroup $L$ of $B$ such that the quotient $B/L$ is
an abelian variety.  By \cite[XVII.7.2.1]{sga3.2} we know that $L$
is a product of a unipotent group $U$ and a torus. But
 there are no non-trivial morphisms of $F$-groups from $U$ to
an abelian variety \cite[2.3]{conrad-chevalley} or to a torus
\cite[XVII.2.4]{sga3.2}, so there are no non-trivial morphisms
from $U$ to $A$. Therefore, $U$ is trivial, and $B$ is
semi-abelian.
\end{proof}
\begin{lemma}\label{lemma-semiab}
If $A$ is a tamely ramified semi-abelian $K$-variety, with N\'eron
model $\mathcal{A}$, then $\widetilde{F}^{>0}\mathcal{A}_s$ is
unipotent, and $(\widetilde{Gr}^0\mathcal{A}_s)^o$ semi-abelian.
\end{lemma}
\begin{proof}
 It follows immediately from Corollary \ref{cor-fixed}  that
$\widetilde{F}^{>0}\mathcal{A}_s$ is unipotent. Let $K'/K$ be a
finite tame extension such that the N\'eron model $\mathcal{A}'$
of $A'=A\times_K K'$ has semi-abelian reduction, and denote by
$\mu$ the Galois group $G(K'/K)$. By Proposition
\ref{prop-tamejump} and Corollary \ref{cor-fixed} we have an
isomorphism
$$\widetilde{Gr}^0\mathcal{A}_s\cong (\mathcal{A}'_s)^{\mu}$$
This is a smooth subgroup scheme of $\mathcal{A}'_s$, and since
$(\mathcal{A}'_s)^o$ is semi-abelian,
$(\widetilde{Gr}^0\mathcal{A}_s)^o$ is semi-abelian by Lemma
\ref{lemma-subsemiab}.
\end{proof}
\begin{cor}\label{cor-lengthunip}
If $A$ is a tamely ramified semi-abelian $K$-variety, with N\'eron
model $\mathcal{A}$, and $(c_1,\ldots,c_v)$ is its tuple of
elementary divisors, then the length $v$ of the tuple equals the
unipotent rank of $\mathcal{A}_s^o$.
\end{cor}
\begin{proof}
By Corollary \ref{cor-compar}, the length $v$ of the tuple of
elementary divisors equals the sum of the multiplicities of the
non-zero jumps of $A$, but this is precisely the dimension of
$\widetilde{F}^{>0}\mathcal{A}_s$.
\end{proof}
%
%

For any integer $i>0$, we denote by $\Phi_i(t)\in \Z[T]$ the
cyclotomic polynomial whose roots are the primitive $i$-th roots
of unity. Its degree equals $\varphi(i)$, with $\varphi(\cdot)$
the Euler function. Recall that we fixed a topological generator
$\sigma$ of the tame monodromy group $G(K^t/K)$, and that we
denote by $\tau:\Q\rightarrow \Z_{>0}$ the function which sends a
rational number to its order in the group $\Q/\Z$.
\begin{theorem}\label{theo-monodromy}
Let $A$ be a tamely ramified abelian $K$-variety, and let
$j_1(A),\ldots,j_u(A)$ be the jumps of $A$, with respective
multiplicities $m_1(A),\ldots,m_u(A)$. Denote by $e$ the degree of
the minimal extension $L$ of $K$ where $A$ acquires semi-abelian
reduction. For each divisor $d$ of $e$ we put
\begin{eqnarray*}
J_d&=&\{i\in \{1,\ldots,u\}\,|\,\tau(j_i(A))=d\}
\\ \nu_d&=&2\cdot
\left( \sum_{i\in J_d} m_i(A)\right)\end{eqnarray*} Then $\nu_d$
is divisible by $\varphi(d)$, and the characteristic polynomial
$P_{\sigma}(t)$ of $\sigma$ on $V_\ell A=T_\ell
A\otimes_{\Z_\ell}\Q_\ell$ is given by
$$P_{\sigma}(t)=\prod_{d|e}\Phi_d(t)^{\nu_d/\varphi(d)}$$
\end{theorem}
\begin{proof}
The polynomial $P_{\sigma}(t)$ belongs to $\Z[t]$, and its zeroes
are roots of unity whose orders divide $e$ \cite[IX.4.3]{sga7a}.
 Hence, $P_{\sigma}(t)$ is a product of cyclotomic polynomials $\Phi_d(t)$ with $d|e$. For
each divisor $d$ of $e$, we denote by $r_d$ the number of zeroes
of $P_{\sigma}(t)$ (counted with multiplicities) that are
primitive $d$-th roots of unity. It suffices to show that
$r_d=\nu_d$ for each divisor $d$ of $e$. We proceed by induction
on $d$.


By Lemma \ref{lemma-semiab}, the multiplicity of $0$ as a jump of
$A$ equals the dimension of the semi-abelian part of
$\mathcal{A}_s^o$. It is well-known that twice this dimension is
equal to the multiplicity of $1$ as an eigenvalue of $\sigma$ on
$V_\ell A$ (see for instance \cite[1.3]{Lenstra-Oort}), so we find
$r_1=\nu_1$.

Now fix a divisor $d$ of $e$ and assume that $\nu_{d'}=r_{d'}$ for
all divisors $d'$ of $d$ with $d'<d$. Note that the multiplicity
$m_d$ of $1$ as a zero of $P_{\sigma^d}(t)$, the characteristic
polynomial of $\sigma^d$ on $V_\ell(A)$, equals
$\sum_{d'|d}r_{d'}$. Applying the previous argument to $A\times_K
K(d)$, we see that $m_d$ equals twice the multiplicity of $0$ as a
jump of $A\times_K K(d)$. By Lemma \ref{lemma-basech}, we obtain
$$\sum_{d'|d}r_{d'}=m_d=\sum_{d'|d}\nu_{d'}$$
so by the induction hypothesis we see that $r_d=\nu_d$.
\end{proof}
\begin{cor}
Let $A$ be a tamely ramified abelian $K$-variety. If $A$ has a
jump $j$ with $\tau(j)=d$, then the primitive $d$-th roots of
unity are monodromy eigenvalues of $\sigma$ on $H^1(A\times_K
K^t,\Q_\ell)$.
\end{cor}

In fact, we can give a more precise description of the relation
between jumps and monodromy. First, we need some auxiliary lemmas.

\begin{lemma}\label{lemma-rep}
For each $n\in \Z_{>0}$ there exists a $\Q$-rational
representation $\rho$ of $\Z/n\Z$ such that $\rho(1)$ has
characteristic polynomial $\Phi_n(t)$.
%
\end{lemma}
\begin{proof}
We proceed by induction on $n$. For $n=1$ the result is clear, so
assume that $n>1$ and that the lemma holds for all $n'<n$. This
implies that, for each divisor $m$ of $n$ with $m<n$, there exists
a $\Q$-rational representation $\rho_{m}$ of $\Z/n\Z$ such that
$\rho_m(1)$ has characteristic polynomial $\Phi_{m}(t)$. Let
$\rho'$ be a complex representation of $\Z/n\Z$ such that
$\rho'(1)$ has characteristic polynomial $\Phi_n(t)$. We'll show
that $\rho'$ is defined over $\Q$. By
\cite[\S\,12.1]{serre-linear}, it is enough to prove that the
character $\chi_{\rho'}$ of $\rho'$ belongs to the representation
ring $R_{\Q}(\Z/n\Z)$. If we denote by $\rho_{reg}$ the
 regular representation of $\Z/n\Z$ over $\Q$, then
$$\rho_{reg}\otimes_{\Q}\C\cong \rho'\oplus
(\bigoplus_{m|n,\,m<n}(\rho_{m}\otimes_{\Q}\C))$$  Hence,
$$\chi_{\rho'}=\chi_{\rho_{reg}}-\sum_{m|n,\,m<n}\chi_{\rho_{m}}$$
belongs to $R_{\Q}(\Z/n\Z)$.
\end{proof}
\begin{lemma}\label{lemma-Z}
Let $A$ be an abelian $K$-variety of dimension $g$, and let
$\gamma$ be an element of $G(K^s/K)$. For any $i\in
\{0,\ldots,2g\}$, the characteristic polynomial
$P_{\gamma}^{(i)}(t)$ of $\gamma$ on $H^i(A\times_K K^s,\Q_\ell)$
belongs to $\Z[t]$. It is independent of $\ell$, and it is a
product of cyclotomic polynomials.
\end{lemma}
\begin{proof}
If $i=1$ then this follows from \cite[IX.4.3]{sga7a}. Hence, by
Lemma \ref{lemma-rep}, there exists an automorphism of finite
order $\gamma'$ of the $\Q$-vector space $V=\Q^{2g}$ such that the
characteristic polynomial of $\gamma'$ on $V$ is
$P_{\gamma}^{(1)}(t)$. For any $i\geq 0$ there exists a
$G(K^s/K)$-equivariant isomorphism of $\Q_\ell$-vector spaces
$$H^i(A\times_K K^s,\Q_\ell)\cong \bigwedge^i H^1(A\times_K
K^s,\Q_\ell)$$ so that $P^{(i)}_{\gamma}(t)$ coincides with the
characteristic polynomial of $\wedge^i \gamma'$ on $\bigwedge^i
V$. It follows that $P_{\gamma}^{(i)}(t)$ is a product of
cyclotomic polynomials, independent of $\ell$.
\end{proof}

\begin{definition}\label{def-tamepol}
We say that a polynomial $Q(t)$ in $\Z[t]$ is $p$-tame if it is
of the form
$$Q(t)=\Phi_{n_1}(t)\cdot \ldots \cdot \Phi_{n_j}(t)$$
with $j\in \Z_{>0}$ and $n_1,\ldots,n_j\in \N'$.
\end{definition}
\begin{lemma}\label{lemm-tamepol}
Consider the unique ring morphism $$\Z[t]\rightarrow k[t]$$
mapping $t$ to $t$.
\begin{enumerate}
 \item Let
$\zeta$ be a primitive $d$-th root of unity in $k$, with $d\in
\N'$. Fix an algebraic closure $\Q^a$ of $\Q$ and a primitive
$d$-th root of unity $\xi$ in $\Q^a$. Let $n$ be an element of
$\Z_{>0}$ and let $a$ be a tuple in $\N^n$. If $Q(t)$ is a
$p$-tame polynomial in $\Z[t]$ whose image in $k[t]$ is divisible
by $\prod_{i=1}^n (t-\zeta^{a_i})$, then $Q(t)$ is divisible by
$\prod_{i=1}^n (t-\xi^{a_i})$ in $\Q^a[t]$.

\item If $Q_1(t)$ and $Q_2(t)$ are $p$-tame polynomials in $\Z[t]$
whose images in $k[t]$ coincide, then $Q_1(t)=Q_2(t)$.
\end{enumerate}
\end{lemma}
\begin{proof}
(1) If $k$ has characteristic zero, the result follows by
considering the embedding of $\Q(\xi)$ in $k$ that maps $\xi$ to
$\zeta$. So assume that $p>1$. Denote by $W(k)$ the ring of
$k$-Witt vectors, and consider the unique embedding of $\Z[\xi]$
in $W(k)$ such that the image of $\xi$ in the residue field $k$ of
$W(k)$ equals $\zeta$. We may assume that the roots of $Q(t)$ are
$d$-th roots of unity. The result then follows from the fact that
the reduction map
$$\mu_d(W(k))\rightarrow \mu_d(k)$$ is a bijection.

(2) Any $p$-tame polynomial $Q(t)$ in $\Z[t]$ divides $(t^u-1)^v$
for some $u\in \N'$ and some $v\in \N$. Hence, the roots of its
image in $k[t]$ are $u$-th roots of unity, so that (2) follows
from (1).
\end{proof}
\begin{prop}\label{prop-abtor}
Let $G$ be a semi-abelian $k$-variety of dimension $g$, and let
$\varphi$ be an automorphism of $G$ of order $d\in \N'$. We fix an
algebraic closure $\Q_\ell^a$ of $\Q_\ell$.
We fix a primitive $d$-th root of unity $\xi$ in $\Q_\ell^a$, and
a primitive $d$-th root of unity $\zeta$ in $k$. Let
$a_1,\ldots,a_g$ be elements of $\{0,\ldots,d-1\}$ such that
$$\zeta^{a_1},\ldots,\zeta^{a_g}$$
are the eigenvalues of $\varphi$ on $Lie(G)$. Then the
characteristic polynomial $Q_{\varphi}(t)$ of $\varphi$ on $V_\ell
G=T_\ell G\otimes_{\Z_\ell}\Q_\ell$ is divisible by
$$\prod_{i=1}^g(t-\xi^{a_i})$$
\end{prop}
\begin{proof}
We consider the Chevalley decomposition
$$0\rightarrow T\rightarrow G\rightarrow B\rightarrow 0$$ of $G$,
with $T$ a $k$-torus and $B$ an abelian $k$-variety. The action of
$\varphi$ on $G$ induces automorphisms on $T$ and $B$, which we
denote again by $\varphi$. There exist $\varphi$-equivariant
isomorphisms
\begin{eqnarray*}
Lie(G)&\cong & Lie(T)\oplus Lie(B)
\\ V_\ell G&\cong &V_\ell T\oplus V_\ell B
\end{eqnarray*}
so that it is enough to consider the case where $G=T$ or $G=B$.

\textit{Case 1: $G=T$.} If we denote by $X(T)$ the character group
of $T$, then there are canonical isomorphisms
\begin{eqnarray}
Lie(T)&=& Hom_{\Z}(X(T),k) \label{eq-lietor}
\\ V_\ell T&=& Hom_{\Z}(X(T),\Q_\ell(1)) \label{eq-tatetor}
\end{eqnarray}
The isomorphism (\ref{eq-tatetor}) implies that $Q_{\varphi}(t)$
belongs to $\Z[t]$. Since $\varphi$ has order $d\in \N'$, we see
that $Q_{\varphi}(t)$ is a $p$-tame polynomial. Combining
isomorphisms (\ref{eq-lietor}) and (\ref{eq-tatetor}), we obtain
that the image of $Q_{\varphi}(t)$ in $k[t]$ equals
$$\prod_{i=1}^g (t-\zeta^{a_i})$$
By Lemma \ref{lemm-tamepol}, we may conclude that
$$Q_{\varphi}(t)=\prod_{i=1}^g (t-\xi^{a_i})$$

\textit{Case 2: $G=B$.} Since $\ell$-adic cohomology is a Weil
cohomology, $Q_{\varphi}(t)$ coincides with the characteristic
polynomial of $\varphi$ on $B$ (see for instance the appendix to
\cite{milne}). In particular, $Q_{\varphi}(t)$ belongs to $\Z[t]$.
Since the order of $\varphi$ belongs to $\N'$, we see that
$Q_{\varphi}(t)$ is a $p$-tame polynomial.

Since $B$ is an abelian variety, its Hodge-de Rham spectral
sequence degenerates at $E_1$ \cite[5.1]{oda}. This yields a
natural short exact sequence
\begin{equation}\label{eq-hodgedr} 0\rightarrow
H^0(B,\Omega^1_B)\rightarrow H^1_{dR}(B)\rightarrow
H^1(B,\mathcal{O}_B)\rightarrow 0\end{equation} with $H^1_{dR}(B)$
the degree one de Rham cohomology of $B$. We have natural
isomorphisms
\begin{equation}\label{eq-lieab}
 H^0(B,\Omega^1_B)\cong
\omega_{B/k}\cong Lie(B)^{\vee}
\end{equation}


First, assume that $k$ has characteristic zero. Then de Rham
cohomology is a Weil cohomology, so that $Q_{\varphi}(t)$ is also
the characteristic polynomial of $\varphi$ on $H^1_{dR}(B)$. By
(\ref{eq-hodgedr}), (\ref{eq-lieab}) and Lemma \ref{lemm-tamepol},
we see that $Q_{\varphi}(t)$ is divisible by
$$\prod_{i=1}^g(t-\xi^{a_i})$$

Now, assume that $k$ has characteristic $p>1$. If we denote by
$H^i_{crys}(B)$ the degree $i$ integral crystalline cohomology of
$B$, then $H^i_{crys}(B)$ is a free $W(k)$-module for each $i$
\cite[II.7.1]{illusie}, so that there is a canonical isomorphism
\begin{equation}\label{eq-crys}
H^1_{dR}(B)\cong H^1_{crys}(B)\otimes_{W(k)}k \end{equation} by
\cite[II.4.9.1]{illusie}. But crystalline cohomology is a Weil
cohomology too, which implies that the characteristic polynomial
of $\varphi$ on $H^1_{crys}(B)$ coincides with  $Q_{\varphi}(t)$.
 By
(\ref{eq-hodgedr}), (\ref{eq-lieab}) and (\ref{eq-crys}), the
image of $Q_{\varphi}(t)$ in $k[t]$ is divisible by
$$\prod_{i=1}^g (t-\zeta^{a_i})$$ By Lemma \ref{lemm-tamepol}, we
see that $Q_{\varphi}(t)$ is divisible by
$$\prod_{i=1}^g (t-\xi^{a_i})$$ in $\Q_\ell^a[t]$.
\end{proof}

\begin{theorem}\label{theo-jumpmon}
Let $A$ be a tamely ramified abelian $K$-variety, and denote by
$e$ the degree of the minimal extension of $K$ where $A$ acquires
semi-abelian reduction. We fix an algebraic closure $\Q^a_\ell$ of
$\Q_\ell$ and a primitive $e$-th root of unity $\xi$ in
$\Q^a_\ell$. If we denote by $j_1(A),\ldots,j_u(A)$ the jumps of
$A$, with respective multiplicities $m_1(A),\ldots,m_u(A)$, then
$$\prod_{i=1}^u (t-\xi^{e\cdot j_i(A)})^{m_i(A)}$$ divides
the characteristic polynomial $P_{\sigma}(t)$ of the action of
$\sigma$ on $V_\ell A=T_\ell A\otimes_{\Z_\ell}\Q_\ell$.
\end{theorem}
\begin{proof}
We put $\mu=G(K(e)/K)$.
 We denote by $(T_\ell A)^{ef}$ the essentially fixed
part of the Tate module $T_\ell A$ \cite[IX.4.1]{sga7a}. It is a
sub-$\Z_\ell$-module of $T_\ell A$, closed under the
$G(K^s/K)$-action. The action of $G(K^s/K)$ on $(T_\ell A)^{ef}$
factors through $\mu$, and we have a canonical $\mu$-equivariant
isomorphism
$$(T_\ell A)^{ef}\cong T_\ell (\mathcal{A}(e)^o_s)$$
by \cite[IX.4.2.6]{sga7a}.

If we denote by $\zeta$ the image of $\sigma$ in $\mu\cong
\mu_e(k)$, then the characteristic polynomial of $\sigma$ on
$Lie(\mathcal{A}(e)_s^o)$ equals
$$\prod_{i=1}^u (t-\zeta^{e\cdot j_i(A)})^{m_i(A)}$$
by Corollary \ref{cor-fixed} and Proposition \ref{prop-tamejump}.
The result now follows from Proposition \ref{prop-abtor}.
\end{proof}

\begin{cor}\label{cor-monocond}
Let $A$ be a tamely ramified abelian $K$-variety of dimension $g$.
The cyclotomic polynomial $\Phi_{\tau(c(A))}(t)$ divides the
characteristic polynomial of $\sigma$ on $H^g(A\times_K
K^t,\Q_\ell)$.
\end{cor}
\begin{proof}
 We keep the notations of Theorem \ref{theo-jumpmon}.
By Corollary \ref{cor-compar} we have
$$c(A)=\sum_{i=1}^{u}m_i(A)\cdot j_i(A)$$
 By Lemma \ref{lemma-Z} it suffices to show that
$$\xi^{e\cdot c(A)}=\prod_{i=1}^u(\xi^{e\cdot
j_i(A)})^{m_i(A)}$$ is an eigenvalue of $\sigma$ on $H^g(A\times_K
K^t,\Q_\ell)$. However, by Theorem \ref{theo-jumpmon}, this is a
product of $\sum_{i=1}^{u}m_i(A)=g$ distinct entries in the
sequence of roots (with multiplicities) of $P_{\sigma}(T)$, and
hence an eigenvalue of $\sigma$ on
$$H^g(A\times_K K^t,\Q_\ell)\cong \bigwedge^g H^1(A\times_K
K^t,\Q_\ell)$$
\end{proof}
\begin{cor}
Let $A$ be a tamely ramified abelian $K$-variety. We use the
notations of Theorem \ref{theo-monodromy}.  For every jump $j$ of
$A$, the multiplicity $m$ of $j$ satisfies $$m\leq \nu_{\tau
(j)}/\varphi(\tau (j))$$
\end{cor}
\begin{proof}
Compare Theorems \ref{theo-monodromy} and \ref{theo-jumpmon}.
\end{proof}

\section{N\'eron models and tame base change}\label{sec-tamebc}
In this section, we assume that $K$ is strictly henselian and that
$k$ is algebraically closed. We adopt the following notation: if
$G$ is a smooth commutative algebraic $K$-group such that $G$
 admits a N\'eron model, then we denote this N\'eron model by
$\mathcal{G}$. If $d\in \N'$ and $G(d)=G\times_K K(d)$ admits a
N\'eron model, we denote it by $\mathcal{G}(d)$.
If $G$ is semi-abelian, we will often use the notations $A$ and
$\mathcal{A}$ instead of $G$ and $\mathcal{G}$.

\begin{definition}\label{def-bc}
Denote by $AV$ the set of isomorphism classes of abelian
$k$-varieties. Let $G$ be a smooth commutative algebraic $K$-group
such that $G(d)$ admits a N\'eron model $\mathcal{G}(d)$ for each
$d\in \N'$. We denote by $u_G(d)$, $t_G(d)$ and $a_G(d)$ the
unipotent, resp. reductive, resp. abelian rank of
$\mathcal{G}(d)_s^o$. We put $\phi_G(d)=\phi(G(d))$ (see
Definition \ref{def-neron}) and we denote by $B_G(d)\in AV$ the
isomorphism class of the abelian quotient in the Chevalley
decomposition of $\mathcal{G}(d)_s^o$.
\end{definition}

\begin{prop}\label{prop-congr}
Let $A$ be a tamely ramified semi-abelian $K$-variety, and let
$j_1(A),\ldots,j_u(A)$ be its jumps, with multiplicities
$m_1(A),\ldots,m_u(A)$. Denote by $e$ the degree of the minimal
extension of $K$ where $A$ acquires semi-abelian reduction. Then
for each $d\in \N'$ we have
$$u_A(d)=\sum_{d\cdot j_i(A)\notin \Z }m_i(A)$$
and this value is also equal to the number of elementary divisors
of $A(d)$. Moreover, the values $u_A(d)$, $t_A(d)$, $a_A(d)$ and
$B_A(d)$ only depend on $d$ mod $e$.
\end{prop}
\begin{proof}
The expressions for $u_A(d)$ follow from Corollaries
\ref{cor-compar} and \ref{cor-lengthunip} and Lemma
\ref{lemma-basech}. Since $e\cdot j_i(A)$ belongs to $\Z$ for each
$i\in \{1,\ldots,u\}$, by Proposition \ref{prop-tamejump}, the
property $d\cdot j_i(A)\in \Z$ only depends on $d\mod e$.

By the equality
$$\mathrm{dim}\,A=u_A(d)+t_A(d)+a_A(d)$$
the only thing left to show is that $B_A(d)\in AV$ only depends on
$d$ mod $e$. By Lemma \ref{lemma-semiab}, it is enough to show
that $(Gr^0\mathcal{A}(d)_s)^o$ and $(Gr^0\mathcal{A}(d+e')_s)^o$
are isomorphic if $d\in \N'$ and $e'$ is a multiple of $e$ such
that $d+e'\in \N'$. We put
\begin{eqnarray*}
m&=&d(d+e')e'
\\ \mu_1&=&G(K(m)/K(d))
\\ \mu_2&=&G(K(m)/K(d+e'))\\ \mu_3&=&G(K(m)/K(e'))\end{eqnarray*}
By Corollary \ref{cor-fixed} we have
\begin{eqnarray*}
(Gr^0\mathcal{A}(d)_s)^o&\cong & ((\mathcal{A}(m)^o_s)^{\mu_1})^o
\\ (Gr^0\mathcal{A}(d+e')_s)^o &\cong & ((\mathcal{A}(m)^o_s)^{\mu_2})^o
\end{eqnarray*}
Let $\zeta$ be a generator of $G(K(m)/K)$. Then $\zeta^d$ and
$\zeta^{d+e'}$ generate $\mu_1$, resp. $\mu_2$. Since $A(e')$ has
semi-abelian reduction, the natural morphism
$\mathcal{A}(e')^o\times_{R(e')}R(m)\rightarrow \mathcal{A}(m)^o$
is an isomorphism. Hence, the action of $\mu_3=<\!\zeta^{e'}\!>$
on $\mathcal{A}(m)^o_s$ is trivial, so the actions of $\zeta^d$
and $\zeta^{d+e'}$ coincide.
\end{proof}

\section{The order function}\label{sec-ord}
%
 In this section, we assume that $K$ is strictly henselian.
 %

\begin{definition}
Let $G$ be a smooth commutative algebraic $K$-group of dimension
$g$, and assume that $G$ admits a N\'eron model $\mathcal{G}$. A
distinguished gauge form on $G$ is a degree $g$ differential form
$\omega$ of the form $\omega=j^*\omega'$, where $j:G\rightarrow
\mathcal{G}$ is the natural open immersion of $G$ into its N\'eron
model, and $\omega'$ is a translation-invariant generator of the
free rank $1$ module $\Omega^g_{\mathcal{G}/R}$.
\end{definition}
Such a distinguished gauge form $\omega$ always exists
\cite[4.2.3]{neron}. It is unique up to multiplication with an
element in the unit group $R^*$, and it is translation-invariant
w.r.t. the group multiplication on $G$.


\begin{definition}
Let $G$ be a smooth commutative algebraic $K$-group, and assume
that $G(d)$ admits a N\'eron model $\mathcal{G}(d)$ for each $d\in
\N'$. Let $\omega$ be a distinguished gauge form on $G$. For each
$d\in \N'$, we put
$$ord_G(d)=-ord(\omega(d))(e_{\mathcal{G}(d)})=-ord_{\mathcal{G}(d)^o_s}(\omega(d))$$
We call $ord_G(\cdot)$ the order function of $G$.
\end{definition}
This definition does not depend on the choice of distinguished
gauge form. The equality
$$ord(\omega(d))(e_{\mathcal{G}(d)})=ord_{\mathcal{G}(d)^o_s}(\omega(d))$$
follows from Proposition \ref{prop-ordsec}. The value $ord_G(d)$
measures the difference between $\mathcal{G}(d)$ and
$\mathcal{G}\times_R R(d)$.
\begin{prop}\label{prop-welldef}
Let $G$ be a smooth commutative algebraic $K$-group, and let $d$
be an element of $\N'$ such that  $G(d)$ admits a N\'eron model
$\mathcal{G}(d)$. If $\omega$ is a gauge form on $G$, then
$$ord_{C}(\omega(d))=ord_{\mathcal{G}(d)_s^o}(\omega(d))$$ for each connected
component $C$ of $\mathcal{G}(d)_s$.
\end{prop}
\begin{proof}
This follows immediately from the translation-invariance of
$\omega$.
\end{proof}
\begin{prop}\label{prop-ord}
Let $G$ be a smooth commutative algebraic $K$-group, and assume
that $G(d)$ admits a N\'eron model $\mathcal{G}(d)$ for each $d\in
\N'$.

 Let $d$ be an element of $\N'$, denote by $h$ the
canonical morphism $\mathcal{G}\times_R R(d)\rightarrow
\mathcal{G}(d)$, and denote by $\mathcal{K}(d)$ the kernel of $h$.
 Then
$$ord_G(d)=
\mathrm{length}_{R(d)}\omega_{\mathcal{K}(d)/R(d)}$$

%
\end{prop}
\begin{proof}
%
By the exact sequence of $R(d)$-modules
$$\begin{CD}0@>>> \omega_{\mathcal{G}(d)/R(d)} @>\alpha>>
\omega_{\mathcal{G}\times_R R(d)/R(d)} @>>>
\omega_{\mathcal{K}(d)/R(d)}@>>> 0\end{CD}$$ from Theorem
\ref{theo-diff}, we have
$$\mathrm{length}_{R(d)}\omega_{\mathcal{K}(d)/R(d)}=\mathrm{length}_{R(d)}coker(\alpha)$$
Since $\alpha$ is an injective morphism of free $R(d)$-modules of
the same rank,
$$\mathrm{length}_{R(d)}coker(\alpha)=\mathrm{length}_{R(d)}coker(det(\alpha))$$
But $det(\alpha)$ is nothing but the morphism of free rank one
$R(d)$-modules
\begin{equation}\label{eq-detalpha}det(\alpha):e_{\mathcal{G}(d)}^*\Omega^g_{\mathcal{G}(d)/R(d)}\rightarrow
e_{\mathcal{G}\times_R R(d)}^*\Omega^g_{\mathcal{G}\times_R
R(d)/R(d)}\end{equation} If we denote by $\omega(d)'$ the unique
extension of $\omega(d)$ to a relative gauge form on
$\mathcal{G}\times_R R(d)$, then the target of (\ref{eq-detalpha})
is generated by the pull-back of $\omega(d)'$. Hence,
$$\mathrm{length}_{R(d)}\omega_{\mathcal{K}(d)/R(d)}=-ord(\omega(d))(e_{\mathcal{G}(d)})=ord_G(d)$$
%
\end{proof}

\begin{prop}\label{prop-computord}
Let $G$ be a smooth commutative algebraic $K$-group, and assume
that $G(d)$ admits a N\'eron model $\mathcal{G}(d)$ for each $d\in
\N'$. If we denote by
$$j_1(G,K(d)),\ldots,j_u(G,K(d))$$ the $K(d)$-jumps of $G$, with
respective multiplicities $m_1(G,K(d)),\ldots,m_u(G,K(d))$, then
we have
$$ord_G(d)=c(G,K(d))\cdot d=\sum_{i=1}^{u}m_i(G,K(d))\cdot j_i(G,K(d))$$
for each $d\in \N'$.

In particular, if $A$ is a tamely ramified semi-abelian
$K$-variety, with jumps $j_1(A),\ldots,j_v(A)$ with respective
multiplicities $m_1(A),\ldots,m_v(A)$, then we have
$$ord_A(d)=\sum_{i=1}^v m_i(A)\cdot \lfloor j_i(A)\cdot d\rfloor $$ for
every $d\in \N'$.
\end{prop}

\begin{proof}
Combine Theorem \ref{theo-diff}, Proposition \ref{prop-tamejump}
and Proposition \ref{prop-ord}.
\end{proof}

\begin{cor}\label{cor-computord}
Let $A$ be a tamely ramified semi-abelian $K$-variety, and denote
by $e$ the degree of the minimal extension of $K$ where $A$
acquires semi-abelian reduction. For all $d\in \N'$ and all $q\in
\Z_{>0}$ such that $d+q\cdot e\in \N'$, we have
$$ord_A(d+q\cdot e)=ord_A(d)+q\cdot c(A)\cdot e$$
Moreover, we have
$$ord_A(d)\leq c(A)\cdot d$$
for all $d\in \N'$, with equality iff $e|d$.
\end{cor}
\begin{proof}
The first assertion follows from the fact that the jumps of $A$
belong to $(1/e)\Z$ (Proposition \ref{prop-tamejump}). For the
second, note that with the notations of Proposition
\ref{prop-computord}, we have
\begin{eqnarray*}
ord_A(d)&=&\sum_{i=1}^v m_i(A)\cdot \lfloor j_i(A)\cdot d\rfloor
\\&\leq& \sum_{i=1}^v m_i(A)\cdot j_i(A)\cdot d
\\ &=&c(A)\cdot d\quad \quad (\mathrm{Proposition}\
\ref{cor-compar})\end{eqnarray*} with equality iff $j_i(A)\cdot
d\in \Z$ for all $i$. By Lemma \ref{lemma-basech}, $j_i(A)\cdot
d\in \Z$ for all $i$ iff all the jumps of $A(d)$ are zero, i.e.,
iff $A(d)$ has semi-abelian reduction (Corollary
\ref{cor-jumps0}), i.e., iff $e|d$.
\end{proof}

%
%

\section{The motivic zeta function of an abelian
variety}\label{sec-rat} Throughout this section, we assume that
$K$ is complete and $k$ algebraically closed. We keep the
notations of Sections \ref{sec-tamebc} and \ref{sec-ord}, in
particular the ones introduced in Definition \ref{def-bc}.
\subsection{The motivic zeta function}
\begin{definition}\label{def-zeta}
Let $G$ be a smooth commutative algebraic $K$-group of dimension
$g$ such that $G$ admits a N\'eron model. Let $\omega$ be a
distinguished gauge form on $G$. We define the motivic zeta
function $Z_G(T)$ of $G$ by
$$Z_G(T)=\LL^g\cdot \sum_{d\in
\N'}\left(\int_{G(d)^{b}}|\omega(d)|\right)T^d \quad \in
\mathcal{M}_k[[T]]$$ where $G(d)^b$ denotes the bounded part of
$G(d)$ (Definition \ref{def-bounded}).

In particular, if $k$ has characteristic zero and $G$ is proper,
 we have
$$Z_G(T)=\LL^g\cdot S(G,\omega;T)\in \mathcal{M}_k[[T]]$$
with $S(G,\omega;T)$ the motivic generating series associated to
$(G,\omega)$ (Section \ref{subsec-motseries}).
\end{definition}
Since $\omega$ is unique up to multiplication with a unit in
$R^*$, this definition is independent of $\omega$. Note that
$G(d)$ admits a N\'eron model $\mathcal{G}(d)$ for each $d\in \N'$
because $R$ is excellent (see Lemma \ref{lemma-sep} and
Proposition \ref{prop-equiv}).

\begin{prop}\label{prop-explicit}
Let $G$ be a smooth commutative algebraic $K$-group of dimension
$g$ such that $G$ admits a N\'eron model $\mathcal{G}$. Let
$\omega$ be a distinguished gauge form on $G$. Then
\begin{eqnarray*}
Z_G(T)&=&\sum_{d\in \N'}[\mathcal{G}(d)_s]\LL^{ord_G(d)}T^d
\\&=&\sum_{d\in \N'}\left(\phi_{G}(d)\cdot (\LL-1)^{t_G(d)}\cdot \LL^{u_G(d)+ord_G(d)}\cdot [B_G(d)]
 \cdot  T^d\right)
\end{eqnarray*}
in $\mathcal{M}_k[[T]]$.
\end{prop}
\begin{proof}
 Since the formal $\pi$-adic completion of $\mathcal{G}(d)$ is a weak
N\'eron model for $G(d)^b$, the first equality follows from
Proposition \ref{prop-motint} and Proposition \ref{prop-welldef}.
The second one follows from \cite[2.1]{Ni-abelian}.
\end{proof}


\subsection{Cohomological interpretation}
\begin{theorem}
Let $A$ be a tamely ramified abelian $K$-variety.  Denote by
$Add_A$ the set of elements $d$ in $\N'$ such that $A(d)$ has
purely additive reduction.

 Then
\begin{eqnarray*}
\chi_{top}(Z_A(T))&=&\sum_{d\in
\N'}\chi_{top}(\mathcal{A}(d)_s)T^d \\&=&\sum_{d\in
Add_A}\phi_A(d) T^d
\\&=&\sum_{d\in \N'}\sum_{i\geq 0}(-1)^i Trace(\sigma^d\,|\,H^i(A\times_K K^t,\Q_\ell))T^d
\end{eqnarray*}
in $\Z[[T]]$.
\end{theorem}
\begin{proof}
The first two equalities follow from Proposition
\ref{prop-explicit}, without assuming that $A$ is tamely ramified.
The last equality follows from the trace formula
\cite[2.5+8]{Ni-abelian}.
\end{proof}
\subsection{Proof of the monodromy conjecture for tamely ramified abelian varieties}

\begin{definition}
Let $A$ be a semi-abelian $K$-variety, and take a finite extension
$K'$ of $K$ such that $A'=A\times_K K'$ has semi-abelian
reduction. We denote by $\mathcal{A}'$ the N\'eron model of $A'$,
and we define the potential toric rank $t_{pot}(A)$ of $A$ to be
the reductive rank of $(\mathcal{A}'_s)^o$. It is independent of
$K'$.
\end{definition}

 The
notion of pole of a rational series in $\mathcal{M}_k[[\LL^{-s}]]$
was defined in \cite[\S\,4]{rove} (here $s$ is a formal variable).
This notion requires some care because $\mathcal{M}_k$ might not
be a domain. The following theorem is the main result of the
present paper.
\begin{theorem}[Monodromy conjecture for abelian varieties]\label{thm-ratzeta}
Let $A$ be a tamely ramified abelian $K$-variety of dimension $g$.
\begin{enumerate}
\item The motivic zeta function $Z_A(T)$ is rational, and belongs
to the subring
$$\mathscr{R}_k^{c(A)}=\mathcal{M}_k\left[T,\frac{1}{1-\LL^a T^b}\right]_{(a,b)\in
\N\times \N_0,\,a/b=c(A)}$$ of $\mathcal{M}_k[[T]]$. The zeta
function $Z_A(\LL^{-s})$ has a unique pole at $s=c(A)$, whose
order is equal to $t_{pot}(A)+1$. The degree of $Z_A(T)$ is equal
to zero if $p=1$ and $A$ has potential good reduction, and
strictly negative in all other cases.

\item The cyclotomic polynomial $\Phi_{\tau(c(A))}(t)$ divides the
characteristic polynomial of the tame monodromy operator $\sigma$
on $H^g(A\times_K K^t,\Q_\ell)$. Hence, for every embedding
$\Q_\ell \hookrightarrow \C$, the value $\exp(2\pi c(A)i)$ is an
eigenvalue of $\sigma$ on $H^g(A\times_K K^t,\Q_\ell)$.
\end{enumerate}
\end{theorem}
\begin{proof}
Statement (2) follows immediately from Corollary
\ref{cor-monocond}, so it suffices to prove (1). In order to
detect poles of $Z_A(T)$, we will specialize $Z_A(T)$ by means of
a ring morphism
$$\rho:\mathcal{M}_k\rightarrow F$$
such that $F$ is a field. We write $c(A)$ as $a_0/b_0$, with
$a_0\in \N$ and $b_0\in \Z_{>0}$. We assume that the image
$\rho(\LL)$ of $\LL$ in $F$ is not a root of unity, and that
$\rho(\LL)$ has a $b_0$-th root $\rho(\LL)_{b_0}$ in $F$. If
$S(T)$ is an element of $\mathscr{R}_k^{c(A)}$ such that
$$\rho(S(T))\in F[[T]]$$ has a pole of order $n$ at
$(\rho(\LL)_{b_0})^{-a_0}$, then it follows easily from the
definition in \cite[\S\,4]{rove} that $S(\LL^{-s})$ has a pole at
$s=c(A)$ of order at least $n$.

The specialization morphism we'll use is the \emph{Poincar\'e
polynomial}
$$P_k:\mathcal{M}_k\rightarrow \Z[u,u^{-1}]$$
(see \cite[\S\,8]{Ni-tracevar}). If $X$ is a $k$-variety of
dimension $d$ with $r$ irreducible components of maximal
dimension, then $P_k([X])$ is a polynomial in $\Z[u]$ of degree
$2d$ whose leading coefficient equals $r$ \cite[8.7]{Ni-tracevar}.
In particular, this leading coefficient is strictly positive. The
element $P_k(\LL)$ is equal to $u^2$. We denote by $\mathscr{F}$
the quotient field of $\cup_{n>0}\Z[u^{1/n},u^{-1/n}]$.

Let $e$ be the degree of the minimal extension of $K$ where $A$
acquires semi-abelian reduction. For each $\alpha\in
\{1,\ldots,e\}$, we put
$$Z^{(\alpha)}_A(T)=\sum_{d\in \N'\cap (\alpha+\N e)}\left(\phi_{A}(d)\cdot (\LL-1)^{t_A(d)}\cdot \LL^{u_A(d)+ord_A(d)}\cdot [B_A(d)]
 \cdot  T^d\right)$$
 By Proposition \ref{prop-explicit}, we have
$$Z_A(T)=\sum_{\alpha=1}^e Z^{(\alpha)}_{A}(T)$$
 Hence, it suffices to prove the following claims.

\textit{Claim 1. For each $\alpha\in \{1,\ldots,e\}$, the series
$Z^{(\alpha)}_A(T)$ belongs to the sub-$\mathcal{M}_k[T]$-module
$\mathcal{M}^{\alpha}$ of $\mathscr{R}^{c(A)}_k$ generated by the
elements
$$(1-\LL^a T^b)^{-t_A(\alpha)-1}$$ with $(a,b)\in \N\times
\N_{0}$, $a/b=c(A)$. Its specialization $P_k(Z^{(\alpha)}_A(T))$
in $\mathscr{F}[[T]]$ has a pole of order $t_A(\alpha)+1$ at
$T=u^{-2 c(A)}$. The residue of this pole belongs to the subring
$\cup_{n>0}\Q[u^{1/n},u^{-1/n}]$ of $\mathscr{F}$, and its leading
coefficient has sign $(-1)^{t_A(\alpha)+1}$.}

\textit{Claim 2. The degree of $Z^{(\alpha)}_A(T)$ is zero if
$p=1$, $\alpha=e$ and $A(e)$ has good reduction, and strictly
negative in all other cases.}

First, we prove Claim 1. For each $\alpha\in \{1,\ldots,e\}$, we
put $\alpha'=gcd(\alpha,e)$. It follows from \cite[5.7]{HaNi-comp}
that
$$\phi_A(d)=(d/\alpha')^{t_A(\alpha')}\phi_A(\alpha')$$
for every element $d$ of $\N'\cap (\alpha+\N e)$. Moreover,
$t_A(\alpha')=t_A(\alpha)$ by \cite[4.2]{HaNi-comp}. Using
Proposition \ref{prop-congr} and Corollary \ref{cor-computord}, we
can write
$$Z^{(\alpha)}_A(T)=\phi_A(\alpha')(\LL-1)^{t_A(\alpha)}\LL^{u_A(\alpha)+ord_A(\alpha)}[B_A(\alpha)]
(\alpha')^{-t_A(\alpha)}S^{(\alpha)}_A(T)$$ with
$$S^{(\alpha)}_A(T)=\sum_{\{q\in \N\,|\,qe+\alpha\in
\N'\}}(qe+\alpha)^{t_A(\alpha)}\LL^{q c(A)e} T^{qe+\alpha}$$ We
denote by $n_\alpha$ the smallest element of $\alpha+\N e$ that is
divisible by $p$, and we put $q_{\alpha}=(n_\alpha-\alpha)/e$.
Note that $n_\alpha \leq pe$, with equality iff $\alpha=e$. We put
$\varepsilon_k=0$ if $p=1$, and $\varepsilon_k=1$ else. With this
notation at hand, we can write $S^{(\alpha)}_A(T)$ as
$$T^{\alpha}\left(\sum_{q\in \N}(qe+\alpha)^{t_A(\alpha)}(\LL^{c(A)} T)^{qe}\right)-
\varepsilon_k \LL^{q_\alpha c(A)e}T^{n_\alpha}\left(\sum_{r\in
\N}(epr+n_\alpha)^{t_A(\alpha)}(\LL^{c(A)}T)^{epr} \right)$$ From
 \cite[6.2]{HaNi-comp} and its proof, we can deduce that
 $S^{(\alpha)}_A(T)$ belongs to $\mathcal{M}^{\alpha}$, and that
 its specialization $$P_k(S^{(\alpha)}_A(T))\in \mathscr{F}[[T]]$$
 has a pole of order $t_A(\alpha)+1$ at $T=u^{-2 c(A)}$, whose
 residue equals
 $$(-1)^{t_A(\alpha)+1}(t_A(\alpha)!)u^{-2c(A)(\alpha+t_A(\alpha)+1)}(e^{-1}
 -\varepsilon_k (ep)^{-1})$$
 This concludes the proof of Claim 1.

 Now we prove claim 2. By our expression for
 $S^{(\alpha)}_A(T)$ and \cite[6.2]{HaNi-comp}, the rational
 function
 $S^{(\alpha)}_A(T)$ has degree $<0$ if $t_A(\alpha)>0$, so that
 we may assume that $t_A(\alpha)=0$. In that case, we find
 $$S^{(\alpha)}_A(T)=T^{\alpha}\left(\sum_{q\in \N}(\LL^{c(A)} T)^{qe}\right)-
\varepsilon_k \LL^{q_\alpha c(A)e} T^{n_\alpha}\left(\sum_{r\in
\N}(\LL^{c(A)}T)^{epr} \right)$$ Direct computation shows that the
degree of this rational function is $\leq 0$, with equality iff
$\alpha=e$ and $p=1$. This concludes the proof.
\end{proof}

\subsection{Elliptic curves}
As an example, we can give an explicit formula for the motivic
zeta function of a tamely ramified elliptic $K$-curve $E$, in
terms of the base change conductor $c(E)$. We refer to
\cite[5.4.5]{edix} and \cite[\S\,8]{halle-neron} for a table with
the values of $c(E)$ (equivalently, the unique jump of $E$) for
each of the Kodaira-N\'eron reduction types.
\begin{prop}\label{prop-elliptic}
 Let $E$ be a tamely
ramified elliptic curve over $K$, and denote by $e$ the degree of
the minimal extension of $K$ where $E$ acquires semi-abelian
reduction.  Denote by $J$ the set of integers in $\{1,\ldots,e
p-1\}$ that are prime to $p$ and not divisible by $e$. In order to
get uniformous formulas, we introduce an error factor
$\varepsilon_k$ which equals zero for $p=1$ and which equals one
for $p>1$. We fix an algebraic closure $\Q^a$ of $\Q$, and denote
by $\xi_1$ and $\xi_2$ the roots in $\Q^a$ of the characteristic
polynomial $P_\sigma(t)\in \Z[t]$ of $\sigma$ on $T_\ell A$.

Then $e=\tau(c(E))$. If $c(E)=0$ then $P_{\sigma}(t)=(t-1)^2$. If
$c(E)=1/2$ then $P_\sigma(t)=(t+1)^2$. If $c(E)\notin \{0,1/2\}$
then $P_{\sigma}(t)=\Phi_{\tau(c(E))}(t)$.

Moreover, if we put
$$S_E(T)=\sum_{i\in J}(1-(\xi_1)^i)(1-(\xi_2)^i)\frac{\LL^{1+\lfloor c(E)i
\rfloor}T^i}{1-\LL^{c(E)ep}T^{ep}}\quad \in \mathcal{M}_k[[T]]$$
then
$$Z_E(T)=[B_E(e)]\cdot \left(\frac{\LL^{c(E)e}T^e}{1-\LL^{c(E)e}T^e}-\varepsilon_k \frac{\LL^{c(E)ep}T^{ep}}{1-\LL^{c(E)ep}T^{ep}}\right)
+S_E(T)$$ if $E$ has potential good reduction, and
$$Z_E(T)=\phi_{E}(e)(\LL-1)\left(\frac{\LL^{c(E)e}T^e}{(1-\LL^{c(E)e}T^e)^2}-\varepsilon_k\frac{p\LL^{c(E)ep}T^{ep}}{(1-\LL^{c(E)ep}T^{ep})^2}\right)
+S_E(T) $$ else.
\end{prop}
\begin{proof}
The equality $e=\tau(c(E))$ and the expressions for
$P_{\sigma}(t)$ follow immediately from Theorem
\ref{theo-monodromy}. If $d$ is an element of $\N'$ that is not
divisible by $e$, then $E(d)$ has additive reduction, so that
$$\phi_E(d)=(1-(\xi_1)^d)(1-(\xi_2)^d)$$ by the trace formula in
\cite[2.8]{Ni-abelian}. This value only depends on the residue
class of $d$ modulo $e$, since $\xi_1$ and $\xi_2$ are $e$-th
roots of unity.

If $d\in \N'$ is divisible by $e$, then $\phi_E(d)=1$ if $E$ has
potential good reduction, and $\phi_E(d)=(d/e)\phi_E(e)$ if $E$
has potential multiplicative reduction \cite[5.7]{HaNi-comp}.
 Now the formulas for $Z_E(T)$ follow easily from
the computation of $ord_E(\cdot)$ in Proposition
\ref{prop-computord}, and the expression for the motivic zeta
function in Proposition \ref{prop-explicit}.
\end{proof}

\section*{Acknowledgements}
The second author is grateful to B. Moonen and W. Veys for
valuable suggestions.


\begin{thebibliography}{10}

\bibitem{sga3.1}
{\em Sch\'emas en groupes. {I}: {P}ropri\'et\'es g\'en\'erales des
sch\'emas en
  groupes}.
\newblock S\'eminaire de G\'eom\'etrie Alg\'ebrique du Bois Marie 1962/64 (SGA
  3). Dirig\'e par M. Demazure et A. Grothendieck. Lecture Notes in
  Mathematics, Vol. 151. Springer-Verlag, Berlin, 1970.

\bibitem{sga3.2}
{\em Sch\'emas en groupes. {II}: Groupes de type multiplicatif, et
structure
  des sch\'emas en groupes g\'en\'eraux}.
\newblock S\'eminaire de G\'eom\'etrie Alg\'ebrique du Bois Marie 1962/64 (SGA
  3). Dirig\'e par M. Demazure et A. Grothendieck. Lecture Notes in
  Mathematics, Vol. 152. Springer-Verlag, Berlin, 1970.

\bibitem{sga7a}
{\em Groupes de monodromie en g\'eom\'etrie alg\'ebrique. {I}}.
\newblock Springer-Verlag, Berlin, 1972.
\newblock S\'eminaire de G\'eom\'etrie Alg\'ebrique du Bois-Marie 1967--1969
  (SGA 7 {I}), Dirig\'e par A. Grothendieck. Avec la collaboration de M.
  Raynaud et D. S. Rim, Lecture Notes in Mathematics, Vol. 288.



\bibitem{Berk-vanish}
V.~G. Berkovich.
\newblock {Vanishing cycles for formal schemes}.
\newblock {\em Invent. Math.}, 115(3):539--571, 1994.

\bibitem{berk-vanish2}
V.~G. Berkovich.
\newblock {Vanishing cycles for formal schemes, {II}}.
\newblock {\em Invent. Math.}, 125(2):367--390, 1996.



\bibitem{neron}
S.~Bosch, W.~{L\"u}tkebohmert, and M.~Raynaud.
\newblock {\em {N\'eron models}}, volume~21 of
 {\em Ergebnisse der Mathematik und ihrer Grenzgebiete}. \newblock Springer-Verlag, 1990.

\bibitem{formner}
S.~Bosch and K.~Schl{\"o}ter.
\newblock N\'eron models in the setting of formal and rigid geometry.
\newblock {\em Math. Ann.}, 301(2):339--362, 1995.

\bibitem{chai}
C.L. Chai.
\newblock {N\'eron models for semiabelian varieties: congruence and change of
  base field}.
\newblock {\em Asian J. Math.}, 4(4):715--736, 2000.

\bibitem{conrad-chevalley}
B.~Conrad.
\newblock {A modern proof of Chevalley's theorem on algebraic groups.}
\newblock {\em J. Ramanujan Math. Soc.}, 17(1):1--18, 2002.


\bibitem{DenefBour}
J.~Denef.
\newblock Report on {I}gusa's local zeta function.
\newblock In {\em S\'eminaire Bourbaki, Vol. 1990/91, Exp. No.730-744}, volume
  201-203, pages 359--386, 1991.

\bibitem{DL5}
J.~Denef and F.~Loeser.
\newblock Motivic {I}gusa zeta functions.
\newblock {\em J. Algebraic Geom.}, 7:505--537, 1998.

\bibitem{DL3}
J.~Denef and F.~Loeser.
\newblock Geometry on arc spaces of algebraic varieties.
\newblock {\em Progr. Math.}, 201:327--348, 2001.

\bibitem{edix}
B.~Edixhoven.
\newblock N\'eron models and tame ramification.
\newblock {\em Compos. Math.}, 81:291--306, 1992.


\bibitem{greenbergII}
M.J.~Greenberg.
\newblock {Schemata over local rings II}.
\newblock {\em Ann. Math.}, 78(2):256--266, 1963.


\bibitem{halle-neron}
L.~Halvard Halle.
\newblock{Galois actions on Neron models of Jacobians.}
\newblock{\em to appear in
Ann. Inst. Fourier}, arXiv:0805.3080.


\bibitem{HaNi-comp}
L.~Halvard Halle and J.~Nicaise.
\newblock{The N\'eron component series of an abelian variety.}
\newblock{\em preprint}, arXiv:0910.1816.

\bibitem{bitt-abelian}
F.~Heinloth.
\newblock {A note of functional equations for zeta functions with values in
  Chow motives.}
\newblock {\em Ann. Inst. Fourier}, 57(6):1927--1945, 2007.

\bibitem{illusie}
L.~Illusie.
\newblock Complexe de de~Rham-Witt et cohomologie cristalline.
\newblock {\em Ann. Sci. {\'E}cole Norm. Sup. (4)}, 12(4):501--661, 1979.






\bibitem{Lenstra-Oort}
H.W. Lenstra and F.~Oort.
\newblock {Abelian varieties having purely additive reduction.}
\newblock {\em J. Pure Appl. Algebra}, 36:281--298, 1985.


\bibitem{liu-lorenzini-raynaud}
Q.~Liu, D.~Lorenzini, and M.~Raynaud.
\newblock {N\'eron models, Lie algebras, and reduction of curves of genus one.}
\newblock {\em Invent. Math.}, 157(3):455--518, 2004.

\bibitem{Loepadic}
F.~Loeser.
\newblock Fonctions d'{I}gusa p-adiques et polyn\^omes de {B}ernstein.
\newblock {\em Am. J. of Math.}, 110:1--22, 1988.

\bibitem{motrigid}
F.~Loeser and J.~Sebag.
\newblock {Motivic integration on smooth rigid varieties and invariants of
  degenerations}.
\newblock {\em Duke Math. J.}, 119:315--344, 2003.






\bibitem{milne}
J.~Milne.
\newblock{Lefschetz classes on abelian varieties.}
\newblock{\em Duke Math. J.}, 96(3), 639--675, 1999.



\bibitem{Ni-tracevar}
J.~Nicaise.
\newblock {A trace formula for varieties over a discretely valued field}.
\newblock {\em to appear in J. Reine Angew. Math.}, arxiv:0805.1323.

\bibitem{Ni-abelian}
J.~Nicaise.
\newblock {Trace formula for component groups of N\'eron models}.
\newblock {\em preprint}, arXiv:0901.1809v2.

\bibitem{ni-trace}
J.~Nicaise.
\newblock A trace formula for rigid varieties, and motivic Weil generating
  series for formal schemes.
\newblock {\em Math. Ann.}, 343(2):285--349, 2009.

\bibitem{Ni-japan}
J.~Nicaise.
\newblock An introduction to $p$-adic and motivic zeta functions and the monodromy conjecture.
\newblock {To appear in the proceedings of the French-Japanese winter school and Zeta and $L$-functions
(Miura, January 2008)}, Memoirs of the Mathematical Society of
Japan, arxiv:0901.4225.

\bibitem{Ni-tori}
J.~Nicaise.
\newblock Motivic invariants of algebraic tori.
\newblock \emph{preprint}, 2009.

\bibitem{NiSe}
J.~Nicaise and J.~Sebag.
\newblock The motivic {S}erre invariant, ramification, and the analytic
  {M}ilnor fiber.
\newblock {\em Invent. Math.}, 168(1):133--173, 2007.

\bibitem{NiSe3}
J.~Nicaise and J.~Sebag.
\newblock Rigid geometry and the monodromy conjecture.
\newblock In D.~Ch\'eniot et~al., editor, {\em Singularity Theory, Proceedings
  of the 2005 Marseille Singularity School and Conference}, pages 819--836.
  World Scientific, 2007.

\bibitem{NiSe-weilres}
J.~Nicaise and J.~Sebag.
\newblock Motivic {S}erre invariants and {W}eil restriction.
\newblock {\em J. Algebra}, 319(4):1585--1610, 2008.

\bibitem{oda}
T.~Oda.
\newblock The first de Rham cohomology group and Dieudonn{\'e} modules.
\newblock {\em Ann. Sci. {\'E}cole Norm. Sup. (4)}, 2(1):63-135, 1969.




\bibitem{Rod}
B.~Rodrigues.
\newblock {On the monodromy conjecture for curves on normal surfaces.}
\newblock {\em Math. Proc. Camb. Philos. Soc.}, 136(2):313--324, 2004.

\bibitem{rove}
B.~Rodrigues and W.~Veys.
\newblock Poles of {Z}eta functions on normal surfaces.
\newblock {\em Proc. London Math. Soc.}, 87(3):164--196, 2003.


\bibitem{sebag1}
J.~Sebag.
\newblock Int\'egration motivique sur les sch\'emas formels.
\newblock {\em Bull. Soc. Math. France}, 132(1):1--54, 2004.

\bibitem{serre-linear}
J.-P. Serre.
\newblock {\em {Repr\'esentations lin\'eaires des groupes finis}}.
\newblock {Paris: Hermann \& Cie}, 1967.

\bibitem{temkin-resol}
M.~Temkin.
\newblock {Desingularization of quasi-excellent schemes in characteristic
  zero.}
\newblock {\em Adv. Math.}, 219(2):488--522, 2008.


\end{thebibliography}
\end{document}